\RequirePackage{etoolbox}
\csdef{input@path}{{style/}{graphics/}}
\documentclass[aap,MSNbibl,dvips]{arximspdf}
\makeatletter
   \@ifpackageloaded{graphicx}{}{\usepackage{graphicx}}
\makeatother

\doi{10.1214/14-AAP1019} 
\volume{25}
\issue{3}
\pubyear{2015}
\firstpage{1155}
\lastpage{1188}

\makeatletter
\newcommand{\rrVert}{\Vert}
\newcommand{\rrvert}{\vert}
\newcommand{\llVert}{\Vert}
\newcommand{\llvert}{\vert}

\newtheorem{teo}{Theorem}
\newtheorem{prop}[teo]{Proposition}
\newproclaim{defn}{Definition}
\newproclaim{asmp}{Assumption}
\newproclaim{exm}{Example}
\makeatother

\begin{document}
\begin{frontmatter}

\title{Multidimensional sticky Brownian motions as limits of exclusion
processes}
\runtitle{Sticky Brownian motions}

\begin{aug}
\author[A]{\fnms{Mikl\'os Z.}~\snm{R\'acz}\thanksref{t1}\ead[label=e1]{racz@stat.berkeley.edu}\ead[label=u1,url]{http://www.stat.berkeley.edu/\textasciitilde racz/}}
\and
\author[A]{\fnms{Mykhaylo}~\snm{Shkolnikov}\corref{}\ead[label=e2]{mshkolni@gmail.com}\ead[label=u2,url]{http://www.stat.berkeley.edu/\textasciitilde mshkolni/}}
\runauthor{M. Z. R\'acz and M. Shkolnikov}
\affiliation{University of California, Berkeley}
\address[A]{Department of Statistics\\
University of California, Berkeley\\
367 Evans Hall\\
Berkeley, California 94720\\
USA\\
\printead{e1}\\
\phantom{E-mail: }\printead*{e2}\\
\printead{u1}\\
\phantom{URL: }\printead*{u2}} 
\end{aug}
\thankstext{t1}{Supported in part by NSF Grant DMS-11-06999 and DOD ONR
Grant N000141110140.}

\received{\smonth{2} \syear{2013}}
\revised{\smonth{12} \syear{2013}}

%
\begin{abstract}
We study exclusion processes on the integer lattice in which particles
change their velocities due to stickiness. Specifically, whenever two
or more particles occupy adjacent sites, they stick together for an
extended period of time, and the entire particle system is slowed down
until the ``collision'' is resolved. We show that under diffusive
scaling of space and time such processes converge to what one might
refer to as a sticky reflected Brownian motion in the wedge. The latter
behaves as a Brownian motion with constant drift vector and diffusion
matrix in the interior of the wedge, and reflects at the boundary of
the wedge after spending an instant of time there. In particular, this
leads to a natural multidimensional generalization of sticky Brownian
motion on the half-line, which is of interest in both queuing theory
and stochastic portfolio theory. For instance, this can model a market,
which experiences a slowdown due to a major event (such as a court
trial between some of the largest firms in the market) deciding about
the new market leader.
\end{abstract}

%
\begin{keyword}[class=AMS]
\kwd{60H30}
\kwd{60K35}
\kwd{60H10}
\end{keyword}
\begin{keyword}
\kwd{Exclusion processes with speed change}
\kwd{reflected Brownian motion}
\kwd{scaling limits}
\kwd{sticky Brownian motion}
\kwd{stochastic differential equations}
\kwd{stochastic portfolio theory}
\end{keyword}

\end{frontmatter}

\section{Introduction}
\label{sec1}

Stochastic processes with sticky points in the Markov process sense
have been studied for more than half a century.
Sticky Brownian motion on the half-line is the process evolving as a
standard Brownian motion away from zero and reflecting at zero after
spending an instant of time there---as opposed to a reflecting Brownian
motion, which reflects instantaneously.
This process was initially studied by Feller~\cite
{feller1952parabolic,feller1954diffusion}, and It{\^o} and McKean~\cite
{ito1996difussion,ito1963brownian} in a more general context, and was
subsequently analyzed in more detail by several further authors~\cite
{harrison1981sticky,amir1991sticky}.
These papers show that sticky Brownian motion arises as a time change
of a reflecting Brownian motion, and that it describes the scaling
limit of random walks on the natural numbers whose jump rate at zero is
significantly smaller than the jump rates at positive sites.

In stochastic analysis, the stochastic differential equation (SDE)
%
%
\begin{equation}
\label{stickyBM} \mathrm{d}S(t)=\mathbf{1}_{\{S(t)>0\}}\, \mathrm{d}B(t)
+ \eta
\mathbf{1}_{\{S(t)=0\}} \,\mathrm{d}t
\end{equation}
satisfied by sticky Brownian motion has drawn much attention, as it is
an example of an SDE for which weak existence and uniqueness hold, but
strong existence and pathwise uniqueness fail; see \cite{Chi}. In fact,
in \cite{warren1997branching} (see also the survey \cite{ES}) it is
shown that a weak solution to (\ref{stickyBM}) cannot be adapted to a
cozy filtration, that is, a filtration generated by a finite or
infinite-dimensional Brownian motion.

The present study is motivated by the question of how one can define
and analyze multidimensional analogues of (\ref{stickyBM}) and whether
solutions to the corresponding systems of SDEs arise as suitable
scaling limits of interacting particle systems in analogy to the
findings of \cite{harrison1981sticky} in the one-dimensional case. In
\cite{karatzas2012systems}, Section~3, it is shown that a large class of
reflecting Brownian motions in the $n$-dimensional wedge
\[
\mathcal{W}=\bigl\{x = (x_1,\ldots, x_n ) \in
\mathbb{R}^n\dvtx  x_1\leq x_2\leq\cdots\leq
x_n\bigr\}
\]
arise as limits of certain exclusion processes with speed change under
diffusive rescaling. As we show below, sticky Brownian motions with
state space $\mathcal{W}$ can also be obtained as scaling limits of
suitable exclusion processes with speed change.

\subsection{Exclusion processes with sticky particles}

To simplify the exposition, we next describe a simple class of particle
systems which converge to sticky Brownian motions in $\mathcal{W}$ in
the scaling limit, and postpone the description of the much wider class
of particle systems that we can handle to Section~\ref{secconv}. We
fix the number of particles $n \in\mathbb{N}$, and also rate
parameters $a >
0$, $\Theta^L = ( \theta_{i,j}^L )_{i\in[n], j\in
[n-1]} \in
[0, \infty)^{n \times( n - 1 )}$, and $\Theta^R = (
\theta_{i,j}^R )_{i\in[n], j\in[n-1]} \in[0, \infty)^{n
\times
( n - 1 )}$ with the notation
\[
[n]=\{1,2,\ldots,n\}.
\]
For a fixed value of the scaling parameter $M > 0$, the particles move
on the rescaled lattice $\mathbb{Z}/\sqrt{M}$; to describe their
motion we
introduce the following Poisson processes, all of which are
independent, and all of which have jump size $\frac{1}{\sqrt{M}}$. For
$i\in[n]$, the Poisson processes $P_i$ and $Q_i$ have jump rates $Ma$,
while for $i\in[n], j\in[n-1]$, the Poisson processes $L_{i,j}$ and
$R_{i,j}$ have jump rates $\sqrt{M} \theta_{i,j}^L$ and $\sqrt{M}
\theta
_{i,j}^R$, respectively. In addition, for notational convenience we
introduce ghost particles at $\pm\infty$, namely: $X_0^M (
\cdot
) \equiv- \infty$ and $X_{n+1}^M ( \cdot) \equiv
\infty$. For any initial condition $X_1^M ( 0 ) < X_2^M
( 0 ) < \cdots< X_n^M ( 0 )$ on $\mathbb{Z}/\sqrt{M}$, we
can then define a particle system evolving on $\mathbb{Z}/\sqrt{M}$ in
continuous time by setting
%
%
\begin{eqnarray}
\label{eqpartsys} \mathrm{d} X_i^M ( t ) & =&
\mathbf{1}_{ \{ X_k^M
( t ) + ({1}/{\sqrt{M}}) < X_{k+1}^M ( t ),
k\in
[n-1] \}} \,\mathrm{d} \bigl( P_i ( t ) -
Q_i ( t ) \bigr)\nonumber
\\
&&{} + \sum_{j=1}^{n-1}
\mathbf{1}_{ \{ X_i^M ( t ) + ({1}/{\sqrt{M}}) < X_{i+1}^M ( t ), X_j^M ( t
)
+ ({1}/{\sqrt{M}}) = X_{j+1}^M ( t ) \}} \,\mathrm{d} R_{i,j} ( t )
\\
&&{} - \sum_{j=1}^{n-1}
\mathbf{1}_{ \{ X_{i-1}^M ( t
)
+ ({1}/{\sqrt{M}}) < X_i^M ( t ), X_j^M ( t
) +
({1}/{\sqrt{M}}) = X_{j+1}^M ( t ) \}} \,\mathrm{d} L_{i,j} ( t ),\nonumber
\end{eqnarray}
for\vspace*{1pt} $i\in[n]$. Note that~(\ref{eqpartsys}) guarantees that for any
$t\geq0$, the particle configuration $(X^M_1(t),X^M_2(t),\ldots,X^M_n(t))$ is an element of the discrete wedge
\[
\mathcal{W}^M= \biggl\{x \in(\mathbb{Z}/\sqrt{M} )^n\dvtx
x_k+\frac
{1}{\sqrt{M}}\leq x_{k+1}, k\in[n-1] \biggr\}.
\]

Intuitively, when apart, the particles move independently on the
rescaled lattice $\mathbb{Z}/\sqrt{M}$ according to the processes $P_i-Q_i$,
$i=1,2,\ldots,n$ (in particular, with jump rates of order $M$);
however, when two particles land on adjacent sites---an event we
describe as a ``collision''---the system experiences a \textit
{slowdown}: the particles change their jump rates to the ones of the
processes $L_{i,j}$ and $R_{i,j}$, $i\in[n], j\in[n-1]$, which are of
order $\sqrt{M}$. The interaction between adjacent particles can be
described as \textit{stickiness}, as it takes a long time (on the time
scale $Mt$) until the collision is resolved and the particles return to
jump rates of order~$M$.

\subsection{Convergence to multidimensional sticky Brownian motions}\label{secintroconv}

The described particle systems converge to a sticky Brownian motion in
$\mathcal{W}$ under the following assumption.
Define $V=(v_{i,j})_{i\in[n], j\in[n-1]}$, the \textit{speed change
matrix}, by setting
%
%
\begin{equation}
\label{eqvij} v_{i,j}:= \cases{ \theta_{i,j}^R -
\theta_{i,j}^L, &\quad if $j \neq i-1,i$,
\vspace*{5pt}\cr
\theta_{i,i-1}^R, &\quad if $j = i-1$,
\vspace*{5pt}\cr
-
\theta_{i,i}^L, &\quad if $j = i$,}
\end{equation}
and\vspace*{1pt} the \emph{reflection matrix} $Q=(q_{j,j'})_{j,j'\in[n-1]}$ by
setting $q_{j,j'} = v_{j+1,j'} - v_{j,j'}$.
When there is a collision between particles $j'$ and $j'+1$ and no
other collisions, then the velocity of particle $i$ is given by
$v_{i,j'}$, and the velocity of gap $j$ between particles $j$ and $j+1$
is given by $q_{j,j'}$.
Define also the $ ( n - 1 ) \times( n - 1 )^2$
matrix $Q^{ ( 2 )} = ( q^{ (2 )}_{i,
(k,\ell)} )_{i,k,l = 1}^{n-1}$ according to
\[
q^{ (2 )}_{i, (k,\ell)}:= \cases{ - \theta_{i,\ell}^L,
&\quad if $k = i-1, \ell\neq i-1$,
\vspace*{5pt}\cr
\theta_{i+1,\ell}^L +
\theta_{i,\ell}^R, &\quad if $k = i, \ell\neq i$,
\vspace*{5pt}\cr
-
\theta_{i+1,\ell}^R, &\quad if $k = i + 1, \ell\neq i + 1$,
\vspace*{5pt}\cr
0, &\quad otherwise.}
\]
Let $q_{\cdot, j}$ denote the $j$th column of $Q$, let $q_{\cdot,
( k, \ell)}^{ (2 )}$ denote the column of $Q^{
(2 )}$ indexed by $ (k,\ell)$ and let $\mathcal
{I}^{
(2 )} \subseteq[n-1 ]^2$ denote the\vspace*{2pt} set of pairs of
indices $ (k,\ell)$ such that $q_{\cdot, ( k, \ell
)}^{ (2 )}$ is the zero vector. Note that $ (k,k )
\in\mathcal{I}^{ (2 )}$ for all $k \in[n-1 ]$.

%
%
\begin{asmp}\label{mainass}
\textup{(a)}
Assume that the matrix $Q$ is \emph{completely-$\mathcal{S}$}, in
the sense that there is a $\lambda\in[0,\infty)^{n-1}$ such that $Q
\lambda\in(0,\infty)^{n-1}$ and the same property is shared by every
principal submatrix of $Q$; see~\cite{TW} for several equivalent definitions.

(b)~Assume that the matrices $Q$ and $Q^{ (2 )}$ (restricted
to nonzero columns) are ``\emph{jointly completely}-$\mathcal{S}$,'' in
the\vspace*{2pt} following sense.
For a vector $u \in\mathbb{R}^k$ and $J \subseteq[k ]$,
let $u^J
\in\mathbb{R}^{\llvert J\rrvert}$ denote the vector obtained from $u$ by
removing all coordinates of $u$ whose index is not in $J$.
We assume that for every $J \subseteq[n-1 ]$, $J \neq
\varnothing$, there exists $\gamma= \gamma( J ) \in(
\mathbb{R}_+ )^{\llvert J\rrvert}$ such that $\gamma\cdot
q_{\cdot, j}^J
\geq1$ for every $j \in J$ and $\gamma\cdot q_{\cdot, (k,\ell
)}^{ (2 ),J} \geq1$ for every $k,\ell\in J$, $
(k,\ell) \notin\mathcal{I}^{ (2 )}$.
\end{asmp}

Under Assumption~\ref{mainass}---which we discuss in more detail
below---we have the following convergence result.

%
%
\begin{teo}\label{mainthm}
Suppose\vspace*{1pt} that Assumption~\ref{mainass} holds, and also that the initial
conditions $ \{ (X^M_1 (0 ),X^M_2 (0
),\ldots,X^M_n (0 ) ), M>0 \}$ are deterministic
and converge to a limit $ (x_1,x_2,\ldots,x_n )\in\mathcal{W}$
as $M\rightarrow\infty$. Then the laws of the paths of the particle
systems $ \{ (X^M_1 (\cdot),X^M_2 (\cdot
),\ldots,X^M_n (\cdot) ), M>0 \}$ on
$D
([0,\infty),\mathbb{R}^n )$ (the space of c\`{a}dl\`{a}g paths with
values in $\mathbb{R}^n$ endowed with the topology of uniform
convergence on
compact sets) converge to the law of the unique weak solution of the
system of~SDEs
%
%
\begin{eqnarray}
\label{mainsde} \qquad \mathrm{d}X_i (t ) &=& \mathbf{1}_{ \{
X_1(t)<X_2(t)<\cdots
<X_n(t) \}}
\sqrt{2a} \,\mathrm{d}W_i (t ) + \sum_{j=1}^{n-1}
\mathbf{1}_{ \{X_j (t )=X_{j+1}
(t
) \}} v_{i,j} \,\mathrm{d}t,
\end{eqnarray}
$i\in[n]$, in $\mathcal{W}$ starting from $(x_1,x_2,\ldots,x_n)$. Here
$ (W_1,W_2,\ldots,W_n )$ is a standard Brownian motion in
$\mathbb{R}^n$.

The solution to (\ref{mainsde}) evolves as a Brownian motion when away
from the boundary $\partial\mathcal{W}$ of $\mathcal{W}$, it does not
spend a nonempty time interval on $\partial\mathcal{W}$; however, it
satisfies
\[
\mathbb{P} \bigl(\mathcal{L} \bigl( \bigl\{t\geq0\dvtx  X(t)\in\partial
{\mathcal
{W}} \bigr\} \bigr)>0 \bigr) = 1,
\]
where $\mathcal{L}$ is the Lebesgue measure on $[
0,\infty)$.
\end{teo}

We refer to the solution of (\ref{mainsde}) with $a=1/2$ as \textit
{sticky Brownian motion} in $\mathcal{W}$ with \textit{reflection
matrix} $V$. We choose this terminology because the SDE~(\ref
{mainsde}) generalizes one-dimensional sticky Brownian motion as in
\cite{harrison1981sticky,amir1991sticky}, and also because it is
consistent with the terminology used in~\cite{TW} and the references
therein dealing with instantaneously reflecting Brownian motions.

Regarding our assumptions, Assumption~\ref{mainass}(a) is a natural
condition, which is necessary for the existence of the limiting
stochastic process; see Theorem~\ref{teouniq}.
Assumption~\ref{mainass}(b) [which is stronger; it implies
Assumption~\ref{mainass}(a)], however, is a technical condition; it is
readily satisfied in many natural situations, but it is not a necessary
condition for the convergence result to hold.
For instance, Assumption~\ref{mainass} is satisfied in the natural
case when $\theta_{i,j}^L = \theta_{i,j}^R = \theta> 0$ for all $i
\in
[n ]$, $j\in[n-1 ]$.
See also~\cite{bhardwaj2009diffusion}, where essentially the same
condition is required and used in the proof of~\cite{bhardwaj2009diffusion}, Theorem~7.7, and where it is shown that this
condition is satisfied if the reflection matrix $Q$ satisfies the
Harrison--Reiman condition~\cite{harrison1981reflected}, and some
additional conditions hold.
On the other hand, consider the case when $\theta_{j,j}^L = \theta
_{j+1,j}^R = \theta> 0$ for all $j \in[n -1 ]$, and
$\theta
_{i,j}^L = \theta_{i,j}^R = 0$ otherwise; in words, suppose that when a
collision occurs, all particles not part of the collision ``freeze,''
that is, they cannot move until the collision is resolved.
It is not hard to see that Assumption~\ref{mainass}(b) cannot hold in
this case, although Assumption~\ref{mainass}(a) holds, and we expect
that the convergence result holds as well.

In Section~\ref{secconv} we prove a much stronger result than
Theorem~\ref{mainthm}, allowing for nonexponential interarrival times
between the jumps in the processes $P_i$, $Q_i$, $L_{i,j}$
and~$R_{i,j}$, as well as for dependence between the latter processes;
see Theorem~\ref{teogen}. This then leads to the definition of a
sticky Brownian motion in $\mathcal{W}$ whose components have unequal
drift and diffusion coefficients. In addition, it is not hard to see
from the proof that for each jump parameter $\theta_{i,j}^L$ or
$\theta
_{i,j}^R$ which is zero, we can choose the jump rate of the
corresponding process $L_{i,j}$ or $R_{i,j}$ to be of order $o
(\sqrt{M} )$ (not necessarily identically zero) for the result of Theorem
\ref{mainthm} to still hold.

One of the main technical difficulties in the proof of Theorem~\ref
{mainthm} and its extension (Theorem~\ref{teogen} in Section~\ref
{secconv}) is posed by the indicator function appearing in the
diffusion matrix of the limiting process. This is in contrast to the
main convergence result in~\cite{karatzas2012systems}, where the
martingale part of the limiting process is a Brownian motion.
Another major difference compared to the setting of~\cite
{karatzas2012systems} is that we consider a large class of
completely-$\mathcal{S}$ reflection matrices and are dealing with weak
solutions of the limiting stochastic differential equation; whereas
in~\cite{karatzas2012systems} only a special class of reflection
matrices is considered, allowing for a pathwise construction of the
limiting object.
Finally, we allow for dependence of interarrival times between jumps
for different particles in Theorem~\ref{teogen} below, which is not
addressed in~\cite{karatzas2012systems}.

\subsection{Applications}

We mention two potential areas of applications for the process in~(\ref
{mainsde}) and its extensions that appear in Section~\ref{seclo}. It
is known that reflected Brownian motions in $\mathcal{W}$ give a class
of tractable descriptive models for the logarithmic market
capitalizations (i.e., the logarithms of the total market values of
stocks) of firms in a large equity market; see, for example, \cite
{karatzas2012systems}. These models lead to realistic capital
distribution curves in the long-run and are also able to produce a
realistic pattern of collisions.
In the same spirit, one can think of~(\ref{mainsde}) as a model for
the logarithmic market capitalizations in an equity market in which the
market experiences a slowdown whenever there is a possibility that two
firms will exchange their ranks (described by a collision).
For example, one can imagine a court trial between two firms, the
result of which decides which firm becomes the market leader, leading
to a slowdown of the market right before the time of the verdict as the
market participants await the result of the trial.
The question of whether real-world equity markets spend a positive
amount of time in collisions [so that the logarithmic market
capitalizations should be modeled by the solution of~(\ref{mainsde})]
or the set of times spent in collisions has zero Lebesgue measure (so
that a reflecting Brownian motion in $\mathcal{W}$ is a more
appropriate model) is a challenging statistical problem which should be
addressed in future research.

Another area of application is the study of diffusion approximations of
storage and queueing networks.
It is well known (see, e.g., the survey~\cite{Wi} and the references
therein) that reflected Brownian motions in the orthant describe the
heavy traffic limits of many queuing networks such as open queuing
networks, single class networks and feedforward multiclass networks.
Moreover, Welch~\cite{Wel} discusses a situation where a customer of a
single server queueing network receives \textit{exceptional service}
when the server is idle before his arrival and standard service when
the server is busy prior to his arrival.
The results in~\cite{Wel}, as well as their extensions to more general
exceptional service policies in~\cite{Lem2} and~\cite{Le} show that the
heavy traffic limits of such networks are described by sticky Brownian
motions on the half-line.
For further information on queuing networks with exceptional service
mechanisms we refer the reader to~\cite{Yeo,harrison1981sticky}
and~\cite{Le}.
Similarly, in the setting of a multi-server queuing network, one can
think of a situation where the servers provide exceptional service to a
customer if the server was idle prior to his arrival, and where such
exceptional service slows down the entire queuing network, for example,
due to a commonly used resource.
In view of the aforementioned results in the single-server case, we
expect \textit{sticky Brownian motions in the orthant}~$ (
\mathbb{R}_+
)^{n-1}$, given by the spacings processes
\[
\bigl(X_2 (\cdot)-X_1 (\cdot),X_3 (\cdot
)-X_2 (\cdot),\ldots,X_n (\cdot)-X_{n-1} (
\cdot) \bigr),
\]
to arise as heavy-traffic limits of multi-server queueing networks with
appropriate exceptional service policies.
We also anticipate the tools developed in this paper to appear at the
heart of the proofs of the corresponding heavy-traffic limit theorems.
In the case that the exceptional service by one of the servers does not
affect other servers, we expect the heavy-traffic limit to be given by
a sticky Brownian motion with a \emph{local} rather than a \emph
{global} slowdown; see Section~\ref{future} for further discussion.

\subsection{Future directions} \label{future}

A natural direction for future work is to study other types of sticky
interaction between particles. Even in the class of exclusion processes
in one dimension, there are avenues to be explored. For instance, the
exclusion processes described by~(\ref{eqpartsys}) experience a
\emph{global} slowdown when a collision occurs, whereas for some
applications it would be interesting to consider particle systems with
\emph{local} slowdown. We believe that the techniques we develop in
Section~\ref{secconv} would carry over to such a setting with
appropriate modifications; however, the difficulty of proving
convergence of such processes to the appropriate continuous object
comes from proving uniqueness for the limiting SDE. We expect the
solution of this SDE to spend a positive amount of time on
lower-dimensional faces of the wedge $\mathcal{W}$, making the analysis
of the process more difficult.

\subsection{Outline}

The rest of the paper is structured as follows. Section~\ref{seclo} is
devoted to the study of sticky Brownian motions in $\mathcal{W}$. In
Section~\ref{secexistuniq} we give the proof of existence and
uniqueness of the weak solution to a system of SDEs generalizing~(\ref
{mainsde}). Then, in Section~\ref{secMarkov} we show that the solution
is a Markov process and study the invariant distributions of a suitably
normalized version thereof. Subsequently, Section~\ref{secconv} deals
with the convergence of exclusion processes to sticky Brownian motions
in $\mathcal{W}$. First in Section~\ref{secconvbasic} we prove
Theorem~\ref{mainthm}, and then in Section~\ref{secgeneral} we state
and prove our main result, namely a generalized version of Theorem~\ref
{mainthm}, which deals with the convergence of exclusion processes
with nonexponential and possibly dependent jump interarrival times to
sticky Brownian motions in $\mathcal{W}$.

\section{Multidimensional sticky Brownian motions} \label{seclo}

This section is devoted to the study of the system of SDEs
%
%
\begin{eqnarray}\label{stickygen}
\mathrm{d}X_i (t ) &=& \mathbf{1}_{ \{X_1 (t
)<X_2 (t )<\cdots<X_n (t ) \}}
\bigl(b_i \,\mathrm{d}t + \mathrm{d}W_i(t) \bigr)
\nonumber\\[-8pt]\\[-8pt]
&&{} + \sum
_{j=1}^{n-1} \mathbf{1}_{
\{X_j (t )=X_{j+1} (t ) \}}
v_{i,j} \,\mathrm{d}t,\nonumber
\end{eqnarray}
$i\in[n]$, where $b_i$, $i\in[n]$, are real constants,
$W=(W_1,W_2,\ldots,W_n)$ is an \mbox{$n$-}dimensional Brownian motion with
zero drift vector and a strictly positive definite diffusion matrix
$\mathfrak{C}=(\mathfrak{c}_{i,i'})_{i,i'\in[n]}$,
$V=(v_{i,j})_{i\in
[n], j\in[n-1]}$ is a\vspace*{1pt} matrix with real entries and the initial
conditions $X_i ( 0 ) = x_i$, $i \in[n ]$, satisfy
$ ( x_1, x_2,\ldots, x_n ) \in\mathcal{W}$. We note that the
diffusion matrix of the process $X$ is both discontinuous and
degenerate, so neither existence nor uniqueness of a weak solution to
(\ref{stickygen}) can be obtained directly from the classical results
in~\cite{SV} or~\cite{BP}.

\subsection{Existence and uniqueness} \label{secexistuniq}

In this subsection, we show that Assumption~\ref{mainass}(a) is
necessary and sufficient for the existence and uniqueness of a weak
solution to~(\ref{stickygen}).
Furthermore, even under Assumption~\ref{mainass}(a) one cannot expect
a strong solution to exist.
Our proof relies on the classical results of~\cite{TW} on the existence
and uniqueness of semimartingale reflecting Brownian motions in an
orthant; this connection highlights the importance of Assumption~\ref
{mainass}(a).
We first recall the main definition and the main result from~\cite{TW}.

\begin{defn}[(\cite{TW}, Definition~1.1)]
Let $\eta\in\mathbb{R}^d$, let $\Gamma$ be a $d \times d$ nondegenerate
covariance matrix, let $R$ be a $d \times d$ matrix and for $i \in
[d ]$, let $F_i = \{ \widetilde{x} \in( \mathbb
{R}_+ )^d\dvtx  \widetilde{x}_i = 0 \}$.
For $\widetilde{x} \in( \mathbb{R}_+ )^d$, a \emph{semimartingale
reflecting Brownian motion} (\emph{SRBM}) \emph{in the orthant} $( \mathbb{R}_+
)^d$ \emph{associated with the data}\vspace*{2pt} $ ( \eta, \Gamma, R )$ \emph{that
starts from} $\widetilde{x}$ is a continuous, $ ( \mathcal{F}_t
)$-adapted, $d$-dimensional process $\widetilde{Z}$ defined on
some filtered probability space $ ( \Omega, ( \mathcal{F}_t
)_{t \geq0}, P_{\widetilde{x}} )$ such that under
$P_{\widetilde{x}}$,
\[
\widetilde{Z} ( t ) = \widetilde{X} ( t ) + R \widetilde{Y} ( t ) \in(
\mathbb{R}_+ )^d \qquad\mbox{for all } t \geq0,
\]
where:
\begin{longlist}[(iii)]
\item$\widetilde{X}$ is a $d$-dimensional Brownian motion with drift
vector $\eta$ and covariance matrix $\Gamma$ such that $ \{
\widetilde{X} ( t ) - \eta t, \mathcal{F}_t, t \geq0
\}$ is a martingale and $\widetilde{X} ( 0 ) = \widetilde
{x}$ $P_{\widetilde{x}}$-a.s.,
\item$\widetilde{Y}$ is an $ (\mathcal{F}_t )$-adapted,
$d$-dimensional process such that $P_{\widetilde{x}}$-a.s. for each $i
\in[d ]$, the $i$th component $\widetilde{Y}_i$ of
$\widetilde{Y}$ satisfies:
\begin{enumerate}[(a)]
\item[(a)] $\widetilde{Y}_i ( 0 ) = 0$,\vspace*{2pt}
\item[(b)] $\widetilde{Y}_i$ is continuous and nondecreasing,\vspace*{2pt}
\item[(c)] $\widetilde{Y}_i$ can increase only when $\widetilde{Z}$ is on
the face $F_i$, that is,
\[
\int_0^t \mathbf{1}_{ \{ \widetilde{Z} ( s ) \in
(
\mathbb{R}_+ )^d \setminus F_i \}} \,\mathrm
{d}\widetilde{Y}_i ( s ) = 0
\]
for all $t \geq0$.
\end{enumerate}
\end{longlist}
$\widetilde{Y}$ is referred to as the ``pushing'' process of
$\widetilde{Z}$.
\end{defn}

%
%
\begin{teo}[(\cite{TW}, Theorem~1.3 and Corollary~1.4)]\label{teoTW}
There\vspace*{1pt} exists a SRBM in the orthant $ ( \mathbb{R}_+ )^d$
with data
$ ( \eta, \Gamma, R )$ that starts from $\widetilde{x}
\in
( \mathbb{R}_+ )^d$ if and only if $R$ is
completely-$\mathcal{S}$.
Moreover, when it exists, the joint law of any SRBM, together with its
associated pushing process, is unique.
\end{teo}

We are now ready to prove our result on the system of SDEs~(\ref{stickygen}).

%
%
\begin{teo}\label{teouniq}
Under Assumption~\ref{mainass}\textup{(a)} there exists a unique weak solution
to~(\ref{stickygen}). Moreover, if Assumption~\ref{mainass}\textup{(a)} does
not hold, there is no weak solution to~(\ref{stickygen}).
\end{teo}

\begin{pf} There are two key ideas in the proof.
The first is to consider the process of spacings
\[
\bigl( X_2 ( \cdot) - X_1 ( \cdot), X_3
( \cdot) - X_2 ( \cdot),\ldots, X_n ( \cdot) -
X_{n-1} ( \cdot) \bigr)
\]
and the process $\sum_{i=1}^n X_i ( \cdot)$, which together
determine the process $X ( \cdot)$.
The second idea is to consider an appropriate (and naturally arising)
time change.
\begin{longlist}[\textit{Step} 2.]
\item[\textit{Step} 1.] We start with the proof of weak existence.
First, from Theorem~\ref{teoTW} it follows that there exists a weak
solution on a suitable filtered probability space $ \{\Omega,
(\mathcal{F}_t )_{t\geq0},P \}$ to the following system
of SDEs:
\[
\mathrm{d}\widehat{Z}_i (t )= (b_{i+1}-b_i )
\,\mathrm{d}t + \mathrm{d}B_i (t ) + \sum
_{j=1}^{n-1} q_{i,j} \,\mathrm{d}
\Lambda_j (t ), \qquad i \in[ n-1 ],
\]
with initial conditions $\widehat{Z}_i ( 0 ) = x_{i+1} -
x_i$, $i \in[ n - 1 ]$, where the vector
$B=(B_1,B_2,\ldots,B_{n-1})$ is a Brownian motion with zero drift vector and diffusion
matrix $A=(a_{i,i'})_{i,i'\in[n-1]}$ given by
%
%
\begin{equation}
\label{eqA} a_{i,i'} = \mathfrak{c}_{i,i'}+
\mathfrak{c}_{i+1,i'+1}-\mathfrak{c}_{i,i'+1}-\mathfrak{c}_{i+1,i'},
\end{equation}
the $\Lambda_j (\cdot)$, $j \in[n-1 ]$, are the
semimartingale local times at zero of the processes $\widehat{Z}_j
(\cdot)$, $j\in[n-1 ]$, respectively, and recall that
$q_{i,j} = v_{i+1,j} - v_{i,j}$.
Note that here we have used the fact that the matrix $Q$ is
completely-$\mathcal{S}$; see Assumption~\ref{mainass}(a).
Next, we can find (after extending the underlying probability space if
necessary) a~Brownian motion $\widehat{\beta}= (\widehat{\beta
}_1,\widehat{\beta}_2,\ldots,\widehat{\beta}_n )$ with zero drift
vector and diffusion matrix $\mathfrak{C}$ such that
\[
B_i(\cdot)=\widehat{\beta}_{i+1}(\cdot)-\widehat{
\beta}_i(\cdot),\qquad i\in[n-1].
\]
Therefore, we can define $\widehat{X}=(\widehat{X}_1,\widehat
{X}_2,\ldots,\widehat{X}_n)$ as the unique process satisfying
\begin{eqnarray*}
&\displaystyle\sum_{i=1}^n \widehat{X}_i(t)
= \sum_{i=1}^n x_i + \sum
_{i=1}^n \Biggl(b_i t+
\widehat{\beta}_i(t)+\sum_{j=1}^{n-1}
v_{i,j}\Lambda_j(t) \Biggr),&
\\
&\displaystyle\bigl(\widehat{X}_2 (t )-\widehat{X}_1 (t ),\ldots,
\widehat{X}_n (t )-\widehat{X}_{n-1} (t ) \bigr) = \bigl(
\widehat{Z}_1 (t ),\ldots,\widehat{Z}_{n-1} (t ) \bigr),&
\end{eqnarray*}
for all $t\geq0$. Finally, we let
\begin{eqnarray*}
T(t) &:=&t+\Lambda(t):=t+\sum_{j=1}^{n-1}
\Lambda_j(t),\qquad t\geq0,
\\
\tau(t) &:=&\inf\bigl\{s\geq0\dvtx  T(s)=t\bigr\},\qquad t\geq0,
\end{eqnarray*}
and set $X (\cdot)=\widehat{X} (\tau(\cdot
) )$. Then clearly
%
%
\begin{equation}
\label{eqX} X_i (\cdot)-X_i (0 )= b_i
\tau(\cdot) +\widehat{\beta}_i \bigl(\tau(\cdot) \bigr) +\sum
_{j=1}^{n-1} v_{i,j}
\Lambda_j \bigl(\tau(\cdot) \bigr),\qquad i\in[n].
\end{equation}
Moreover, we note that $\tau(\cdot)$, $\Lambda
(\tau
(\cdot) )$ are nondecreasing functions, which induce
nonnegative measures $\mathrm{d}\tau(\cdot)$, $\mathrm
{d}\Lambda(\tau(\cdot) )$ on $[0,\infty)$
satisfying
%
%
\begin{eqnarray}
\label{measureid} \mathrm{d}\tau(t)+\mathrm{d}\Lambda\bigl(\tau(t)\bigr)=
\mathrm{d}t.
\end{eqnarray}
Therefore we have
\begin{eqnarray*}
\tau(\cdot) &=& \int_0^{\tau(\cdot)} \mathbf{1}_{ \{\widehat
{X}_1
(t )<\widehat{X}_2 (t )<\cdots<\widehat{X}_n
(t
) \}}
\,\mathrm{d}t  = \int_0^\cdot
\mathbf{1}_{ \{X_1 (t )<X_2
(t
)<\cdots<X_n (t ) \}} \,\mathrm{d}\tau(t )
\\
& =& \int_0^\cdot\mathbf{1}_{ \{X_1 (t )<X_2
(t
)<\cdots<X_n (t ) \}}
\,\mathrm{d}t,
\end{eqnarray*}
which takes care of the first term on the right-hand side of~(\ref{eqX}).
In addition, the processes $\widehat{\beta}_i(\tau(\cdot))$, $i\in
[n ]$ are martingales with respect to the filtration $(\mathcal
{F}_{\tau(t )})_{t\geq0}$ with quadratic covariation
processes given by
\[
\mathfrak{c}_{i,i'} \tau(\cdot) = \mathfrak{c}_{i,i'} \int
_0^\cdot\mathbf{1}_{ \{X_1 (t )<X_2 (t )<\cdots
<X_n
(t ) \}} \,\mathrm{d}t,
\qquad i,i'\in[n],
\]
where we used the identity derived in the previous display.
From the last computation we can conclude, in particular, that after
extending the underlying probability space if necessary, we can find a
Brownian motion $\beta= (\beta_1,\beta_2,\ldots,\beta_n )$
with zero drift vector and diffusion matrix $\mathfrak{C}$ such that
\[
\widehat{\beta}_i\bigl(\tau(\cdot)\bigr)=\int_0^\cdot
\mathbf{1}_{ \{
X_1
(t )<X_2 (t )<\cdots<X_n (t ) \}} \,\mathrm{d}\beta_i (t ),\qquad i\in[n].
\]
Finally, we have
\begin{eqnarray*}
\Lambda_j \bigl(\tau(\cdot) \bigr) & =&\int_0^\cdot
\mathbf{1}_{ \{\widehat{X}_j (\tau(t )
)=\widehat
{X}_{j+1} (\tau(t ) ) \}} \,\mathrm{d}\Lambda_j \bigl(\tau(t )
\bigr)
\\
& =&\int_0^\cdot\mathbf{1}_{ \{\widehat{X}_j (\tau
(t ) )=\widehat{X}_{j+1} (\tau(t
) )
\}}
\bigl(\mathrm{d}t-\mathrm{d}\tau(t ) \bigr)
\\
& =&\int_0^\cdot\mathbf{1}_{ \{X_j (t )=X_{j+1}
(t ) \}}
\,\mathrm{d}t - \int_0^{\tau(\cdot)} \mathbf{1}_{ \{\widehat{X}_j (t
)=\widehat
{X}_{j+1}
(t ) \}}
\,\mathrm{d}t
\\
& =&\int_0^\cdot\mathbf{1}_{ \{X_j (t )=X_{j+1}
(t ) \}}
\,\mathrm{d}t,
\end{eqnarray*}
for $j\in[n-1]$. Here the second identity is a consequence of (\ref
{measureid}) and the fact that the boundary local times $\Lambda_{j'}$,
$j'\neq j$, do not charge the set $\{t\dvtx  \widehat{X}_j(t)=\widehat
{X}_{j+1}(t)\}$ (see the main result, Theorem~1, in~\cite{RW}); and the
fourth identity follows from the fact that the instantaneously
reflecting Brownian motion $\widehat{Z}$ does not spend time on the
boundary of the orthant $ (\mathbb{R}_+ )^{n-1}$ \cite{TW}, Lemma~2.1. All in all, we can now conclude that $(X,\beta)$ is a
weak solution to (\ref{stickygen}).

\item[\textit{Step} 2.] We now turn to the proof of weak uniqueness.
To this end, let $(X,W)$ be any weak solution to (\ref{stickygen}). Define
\[
\sigma(t)=\inf\biggl\{s\geq0\dvtx  \int_0^s
\mathbf{1}_{ \{
X_1
(a )<X_2 (a )<\cdots<X_n (a ) \}} \,\mathrm{d}a = t \biggr\},\qquad t\geq0,
\]
and set $\widehat{X} (\cdot)=X (\sigma(\cdot
) )$. Using L\'{e}vy's characterization of Brownian motion, one
verifies that
\[
\widehat{X}_i (t )=\widehat{X}_i (0 )+b_i
t+\widehat{W}_i (t )+\sum_{j=1}^{n-1}
v_{i,j} L_j (t ),\qquad t\geq0,
\]
where $\widehat{W}=(\widehat{W}_1,\widehat{W}_2,\ldots,\widehat{W}_n)$
is a Brownian motion with zero drift vector and diffusion matrix
$\mathfrak{C}$, and $ \{ L_j \}_{j\in[n-1 ]}$ are
nondecreasing processes whose points of increase are contained in the sets
\[
\bigl\{t\geq0\dvtx  \widehat{X}_j (t )=\widehat{X}_{j+1} (t )
\bigr\},\qquad j\in[n-1],
\]
respectively. Moreover, the law of $\widehat{X}$ is uniquely determined
by the joint law of
%
%
\begin{equation}
\label{2proc} \qquad\bigl(\widehat{X}_2 (\cdot)-\widehat{X}_1
(\cdot),\widehat{X}_3 (\cdot)-\widehat{X}_2 (\cdot),
\ldots,\widehat{X}_n (\cdot)-\widehat{X}_{n-1} (\cdot)
\bigr) \quad\mbox{and}\quad\sum_{i=1}^n
\widehat{X}_i(\cdot).
\end{equation}
However, by the uniqueness result of Theorem~\ref{teoTW} we can
identify the first of the latter two processes as an instantaneously
reflected Brownian motion in the orthant $ (\mathbb{R}_+
)^{n-1}$, so
the joint law of that process and its boundary local times is uniquely
determined. Moreover, the second process can be constructed by using
the first process, its boundary local time processes and an additional
independent one-dimensional standard Brownian motion, so the joint law
of the processes in (\ref{2proc}) is uniquely determined. Thus, the law
of $\widehat{X}$ is uniquely determined as well. Finally, the law of
$X$ is also uniquely determined as one can verify that $X (\cdot
)=\widehat{X} (\tau(\cdot) )$, where
$\tau$
is defined as in step~1 above.

\item[\textit{Step} 3.] Suppose now that Assumption~\ref{mainass}(a)
does not hold. Then a weak solution of~(\ref{stickygen}) cannot exist.
Indeed, if $ (X,W )$ was such a weak solution, we could define
the time change $\sigma(\cdot)$ as in step~2 above and let $\widehat
{X}(\cdot)=X(\sigma(\cdot))$ as before. Then the arguments in step~2
would show that the process of spacings
\[
\bigl(\widehat{X}_2 (\cdot)-\widehat{X}_1 (\cdot),
\widehat{X}_3 (\cdot)-\widehat{X}_2 (\cdot),\ldots,
\widehat{X}_n (\cdot)-\widehat{X}_{n-1} (\cdot) \bigr)
\]
is a reflecting Brownian motion in the orthant $ (\mathbb
{R}_+
)^{n-1}$ in the sense of~\cite{TW}. However, by Theorem~\ref{teoTW}
the latter process does not exist if the reflection matrix $Q$ is not
completely-$\mathcal{S}$. This is the desired contradiction.\quad\qed
\end{longlist}\noqed
\end{pf}

The following example shows that, even when Assumption~\ref
{mainass}(a) holds, one cannot expect a strong solution to (\ref
{stickygen}) to exist.

%
%
\begin{exm}
Consider\vspace*{1pt} the following specification of parameters: $n=2$, $b_1=b_2=0$,
$\mathfrak{c}_{1,1}=\mathfrak{c}_{2,2}=1$, $\mathfrak
{c}_{1,2}=\mathfrak
{c}_{2,1}=0$, $v_{1,1}=-\frac{1}{2}$, $v_{2,1}=\frac{1}{2}$; that is,
the system of SDEs is
%
%
\begin{eqnarray}
\mathrm{d}X_1(t) & =&\mathbf{1}_{ \{X_1 (t )<X_2
(t ) \}}
\,\mathrm{d}W_1(t) - \tfrac{1}{2} \mathbf{1}_{ \{
X_1 (t )=X_2 (t ) \}}
\,\mathrm{d}t, \label
{X1}
\\
\mathrm{d}X_2(t) & =&\mathbf{1}_{ \{X_1 (t )<X_2
(t ) \}}
\,\mathrm{d}W_2(t) + \tfrac{1}{2} \mathbf{1}_{ \{
X_1 (t )=X_2 (t ) \}}
\,\mathrm{d}t, \label{X2}
\end{eqnarray}
with $W_1$ and $W_2$ being independent one-dimensional standard
Brownian motions. We claim that this system does not admit a strong
solution. It is well known (see Theorem~3.2 in~\cite{Che}) that strong
existence and weak uniqueness together imply pathwise uniqueness, so it
suffices to show that pathwise uniqueness does not hold for the
system~(\ref{X1})--(\ref{X2}). To this end, we consider the SDE
%
%
\begin{equation}
\label{auxeq} \mathrm{d}Z(t)=\mathbf{1}_{ \{Z (t )>0 \}} \,\mathrm
{d}\beta(t ) +
\mathbf{1}_{ \{Z (t
)=0 \}} \,\mathrm{d}t,
\end{equation}
where $\beta$ is a Brownian motion with zero drift and diffusion
coefficient $2$. The main result in~\cite{Chi} shows that pathwise
uniqueness does not hold for this equation. Therefore it suffices to
argue that pathwise uniqueness for the system (\ref{X1})--(\ref{X2})
would imply pathwise uniqueness for equation~(\ref{auxeq}). Indeed, let
$Z$, $Z'$ be two solutions of~(\ref{auxeq}) on the same probability
space and with respect to the same Brownian motion $\beta$. Extend the
probability space so that it supports an independent Brownian motion
$W$ with zero drift and diffusion coefficient $2$, and define $S$, $S'$
according to
\[
\mathrm{d}S (t )=\mathbf{1}_{ \{Z (t
)>0 \}
} \,\mathrm{d}W (t )\quad\mbox{and}\quad
\mathrm{d}S' (t )=\mathbf{1}_{ \{Z' (t )>0 \}} \,\mathrm{d}W(t).
\]
Finally, set
\[
X_1=\frac{S-Z}{2},\qquad X_2=\frac{S+Z}{2}
\]
and
\[
X'_1=\frac
{S'-Z'}{2},\qquad
X'_2=\frac{S'+Z'}{2}.
\]
Then both $ (X_1,X_2 )$ and $ (X'_1,X'_2 )$ are weak
solutions of the system (\ref{X1})--(\ref{X2}) with respect to the
Brownian motion $ ( (W-\beta)/2, (W+\beta
)/2 )$. Therefore if pathwise uniqueness did hold for the system
(\ref{X1})--(\ref{X2}), we would be able to conclude that $X_1=X'_1$
and $X_2=X'_2$ pathwise, and, hence, that $Z=Z'$ pathwise; in other
words, the solution of (\ref{auxeq}) would be pathwise unique. This is
the desired contradiction.
\end{exm}

\subsection{Markov property and invariant measures} \label{secMarkov}

Having established that the weak solution $X$ of the system~(\ref
{stickygen}) exists and is unique (see Theorem~\ref{teouniq}), we can
now proceed to study some of its properties.
First, we remark that weak existence and uniqueness imply that the
corresponding martingale problem is well posed; see, for example,
Corollary~4.8 and Corollary~4.9 in Chapter~5 of~\cite{KS}.
Therefore, by Theorem~6.2.2 in~\cite{SV}, the process $X$ is Markovian.
In addition, the relation $X(\cdot)=\widehat{X}(\tau(\cdot))$, where
$\widehat{X}$ is an instantaneously reflecting Brownian motion in the
wedge $\mathcal{W}$ with a nondegenerate diffusion matrix, shows that
the process $X$ has the Harris property (see, e.g., the Appendix
of~\cite{dupuis1994lyapunov}),
%
%
\begin{eqnarray}
\label{harris} \forall x, y\in\mathcal{W}, r>0:\qquad\mathbb{P}^x
\bigl(\bigl\llvert X (t )-y\bigr\rrvert<r\mbox{ for some }t\geq0 \bigr)>0.
\end{eqnarray}
Moreover, the corresponding property is true for the process of spacings
\[
Z (\cdot)= \bigl(X_2 (\cdot)-X_1 (\cdot
),X_3 (\cdot)-X_2 (\cdot),\ldots,X_n (
\cdot)-X_{n-1} (\cdot) \bigr).
\]
Thus $Z$ has a unique invariant distribution provided that it is
positive recurrent or, equivalently, if
\[
\widehat{Z} (\cdot)= \bigl(\widehat{X}_2 (\cdot)-
\widehat{X}_1 (\cdot),\widehat{X}_3 (\cdot)-
\widehat{X}_2 (\cdot),\ldots,\widehat{X}_n (\cdot)-
\widehat{X}_{n-1} (\cdot) \bigr)
\]
is positive recurrent; see~\cite{meyn1993stability} and the references
therein. By Proposition~2.8 in the dissertation~\cite{I}, the latter is
the case if and only if
%
%
\begin{eqnarray}
\label{recurrent} Q^{-1} (b_2-b_1,b_3-b_2,
\ldots,b_n-b_{n-1} )^T<0
\end{eqnarray}
componentwise. Here, the superscript $T$ stands for the transpose of
the vector under consideration. We summarize our findings in the next
proposition.

%
%
\begin{prop}\label{propMarkov}
The processes $X$ and $Z$ are Markovian. Both of them possess the
Harris property. Moreover, the process $Z$ has a unique invariant
distribution if and only if the recurrence condition (\ref{recurrent})
is satisfied.
\end{prop}

For a wide class of coefficients the invariant distribution of the
process $Z$ can be given explicitly.
Let $ \langle\cdot, \cdot\rangle$ denote the standard inner
product in Euclidean space, that is, for $x,y \in\mathbb{R}^d$,
$ \langle
x,y \rangle= \sum_{i=1}^d x_i y_i$,
and let $F_i:= \{ z \in( \mathbb{R}_+ )^{n-1}\dvtx
z_i = 0
\}$\vspace*{2pt} denote the $i$th face of the orthant $ ( \mathbb{R}_+
)^{n-1}$.
With this notation we then have the following result.

\begin{teo}
Suppose that in addition to~(\ref{recurrent}) the condition
%
%
\begin{equation}
\label{eqinvcond} 2A=QD+DQ^T
\end{equation}
is satisfied, where $A$ is given by~(\ref{eqA}) and $D=\operatorname{diag}(A)$ (the diagonal matrix, whose diagonal elements coincide with
those of $A$).
Let
%
%
\begin{equation}
\label{eqgamma} \gamma=2 D^{-1}Q^{-1} (b_2-b_1,b_3-b_2,
\ldots,b_n-b_{n-1} )^T,
\end{equation}
and write $\gamma= ( \gamma_1,\ldots, \gamma_{n-1} )$.
Then the invariant distribution of the process $Z$ in the orthant
$
( \mathbb{R}_+ )^{n-1}$ is given by
%
%
\begin{equation}
\label{eqinvdist} \frac{1}{C} e^{ \langle\gamma,z \rangle} \Biggl(
\mathrm{d}z + \sum
_{j=1}^{n-1} \frac{\sqrt{a_{j,j}}}{2}
\mathbf{1}_{ \{ z
\in
F_j \}} \,\mathrm{d}z^j \Biggr),
\end{equation}
where $C=\frac{1 - \sum_{j=1}^{n-1} \sqrt{a_{j,j}} \gamma_j /
2}{\prod_{j=1}^{n-1} ( - \gamma_j )}$ is the appropriate
normalization constant,
and $\mathrm{d}z^j$, $j \in[n-1 ]$, are the Lebesgue
boundary measures on the faces $F_j$, $j \in[n-1 ]$, respectively.
\end{teo}

\begin{pf}
We first transform our process of interest $Z$ in such a way as to make
the diffusion matrix of the transformed process the identity; we refer
to~\cite{ichiba2010collisions}, Section~3.2.1, for similar computations.
Let $U = ( u_{i,j} )_{i,j \in[n-1 ]}$ be an
orthogonal matrix whose columns are the orthonormal eigenvectors of
$A$, and let $G:= U^T A U$, a diagonal matrix with the eigenvalues of
$A$ in its diagonal. Define the transformed process
\[
\overline{Z} ( \cdot):= G^{-1/2} U^T Z ( \cdot).
\]
This is a sticky Brownian motion in the cone
\[
\overline{\mathcal{S}}:= G^{-1/2} U^T ( \mathbb{R}_+
)^{n-1} = \bigl\{ \overline{z} \in\mathbb{R}^{n-1}\dvtx  U
G^{1/2} \overline{z} \in( \mathbb{R}_+ )^{n-1} \bigr\},
\]
with drift vector $\overline{\mu}:= G^{-1/2} U^T \mu$, where $\mu:=
( b_2 - b_1,\ldots, b_n - b_{n-1} )^T$, identity diffusion
matrix, and reflection matrix $\overline{Q}:= G^{-1/2} U^T Q$. This
transformed reflection matrix can be decomposed as
\[
\overline{Q} = ( \overline{\mathfrak{N}} + \overline{\mathfrak{T}} )
D^{-1/2} \equiv( \overline{q}_{\cdot, 1},\ldots, \overline
{q}_{\cdot, n-1} ),
\]
where
\begin{eqnarray*}
\overline{\mathfrak{N}} &:=& G^{1/2} U^T D^{-1/2}
\equiv( \overline{\mathfrak{n}}_{\cdot, 1},\ldots, \overline{\mathfrak
{n}}_{\cdot, n-1} ),
\\
\overline{\mathfrak{T}} &:=& G^{-1/2} U^T Q
D^{1/2} - \overline{\mathfrak{N}} \equiv( \overline{
\mathfrak{t}}_{\cdot, 1},\ldots, \overline{\mathfrak{t}}_{\cdot, n-1} ).
\end{eqnarray*}
The columns of $\overline{\mathfrak{N}}$ are unit vectors, since for
all $j \in[n-1 ]$ we have
%
%
\begin{equation}
\label{eqnunit} \sum_{i=1}^{n-1} \overline{
\mathfrak{n}}_{i,j}^2 = \sum_{i=1}^{n-1}
\biggl( \sqrt{g_{i,i}} u_{j,i} \frac{1}{\sqrt{a_{j,j}}}
\biggr)^2 = \frac
{1}{a_{j,j}} \sum_{i=1}^{n-1}
g_{i,i} u_{j,i}^2 = 1,
\end{equation}
where $G = ( g_{i,j} )_{i,j \in[n-1 ]}$, and in
the last equality we used that $U G U^T = A$. Furthermore, the
corresponding columns of $\overline{\mathfrak{N}}$ and $\overline
{\mathfrak{T}}$ are orthogonal, since for every $i \in
[n-1
]$ we have
\begin{eqnarray*}
\overline{\mathfrak{n}}_{\cdot, i}^T \overline{\mathfrak
{t}}_{\cdot,
i} & =& \sum_{j=1}^{n-1}
\overline{\mathfrak{n}}_{j,i} \overline{\mathfrak{t}}_{j,i} =
\sum_{j=1}^{n-1} \sqrt{g_{j,j}}
u_{i,j} \frac
{1}{\sqrt{a_{i,i}}} \sum_{k=1}^{n-1}
\frac{1}{\sqrt{g_{j,j}}} u_{k,j} q_{k,i} \sqrt{a_{i,i}} -
\sum_{j=1}^{n-1} \overline{\mathfrak
{n}}_{j,i}^2
\\
& =& \sum_{k=1}^{n-1} q_{k,i}
\sum_{j = 1}^{n-1} u_{i,j}
u_{k,j} - 1 = \sum_{k=1}^{n-1}
q_{k,i} \mathbf{1}_{ \{ k = i \}} - 1 = q_{i,i} - 1 = 0,
\end{eqnarray*}
where we used~(\ref{eqnunit}), the fact that $U$ is orthogonal, and
that $\operatorname{diag} ( Q ) = I$, which follows
from~(\ref{eqinvcond}). In fact, $\overline{\mathfrak{n}}_{\cdot,
i}$ is the
inward unit normal to the $i$th face $\overline{F}_i:= G^{-1/2} U^T
F_i$ of the new state space $\overline{\mathcal{S}}$. To see this, let
$w \in F_i$, and let $v:= G^{-1/2} U^T w \in\overline{F}_i$. Then
\begin{eqnarray*}
\overline{\mathfrak{n}}_{\cdot, i}^T v & =& \sum
_{j=1}^{n-1} \overline{\mathfrak{n}}_{j,i}
v_j = \sum_{j=1}^{n-1}
\sqrt{g_{j,j}} u_{i,j} \frac
{1}{\sqrt{a_{i,i}}} \sum
_{k=1}^{n-1} \frac{1}{\sqrt{g_{j,j}}} u_{k,j}
w_k
\\
& =& \frac{1}{\sqrt{a_{i,i}}} \sum_{k=1}^{n-1}
w_k \sum_{j=1}^{n-1}
u_{i,j} u_{k,j} = \frac{1}{\sqrt{a_{i,i}}} \sum
_{k=1}^{n-1} w_k \mathbf
{1}_{ \{ k = i \}}
= \frac{1}{\sqrt{a_{i,i}}} w_i = 0,
\end{eqnarray*}
where the last equality is because $w \in F_i$. Thus the $i$th column
$\overline{q}_{\cdot,i}$ of the new reflection matrix $\overline{Q}$ is
decomposed into components that are normal and tangential to $\overline{F}_i$,
%
%
\begin{equation}
\label{eqqdecomp} \overline{q}_{\cdot,i} = \frac{1}{\sqrt{a_{i,i}}} (
\overline{
\mathfrak{n}}_{\cdot, i} + \overline{\mathfrak{t}}_{\cdot,
i} ).
\end{equation}

The advantage of this transformation is that the setup of the new
process $\overline{Z}$ fits precisely into the framework of Harrison
and Williams~\cite{HW}, who studied the stationary distribution of
reflected Brownian motion with identity diffusion matrix in a convex
polyhedral domain.
Their main result is that the stationary distribution is of exponential
form if the reflection matrix satisfies a certain skew symmetry
condition, and they give explicit formulas for the exponent.
The main difference between their setting and ours is that the process
we study is sticky at the boundary of the domain, as opposed to
reflecting instantaneously, as is the case in~\cite{HW}.
However, apart from taking care of this distinction at the boundary,
the same methods and computations apply.

In particular, we can plug our expressions into the formulas of
Harrison and Williams~\cite{HW} to arrive at the skew symmetry
condition for our process, and also to find the appropriate exponent.
First, we have
%
%
\begin{equation}
\label{eqskew} \overline{\mathfrak{N}}^T \overline{\mathfrak{T}} +
\overline{\mathfrak{T}}^T \overline{\mathfrak{N}} =
D^{-1/2} \bigl( Q D + D Q^T - 2A \bigr) D^{-1/2}.
\end{equation}
The skew symmetry condition of Harrison and Williams (see~\cite{HW}, equation~(1.3)) says that the left-hand side of~(\ref{eqskew})
is the zero matrix, which is the same as our condition~(\ref{eqinvcond}).
Second, plugging in to~\cite{HW}, equation~(4.9), the $\overline
{\gamma
}$ arising in the exponent of the stationary distribution of
exponential form for the transformed process $\overline{Z}$ should be
%
%
\begin{equation}
\label{eqgammaformula} \overline{\gamma} = 2 \bigl( I - \bigl(
\overline{
\mathfrak{N}}^T \bigr)^{-1} \overline{\mathfrak{T}}^T
\bigr)^{-1} \overline{\mu} = 2 G^{1/2} U^T
D^{-1} Q^{-1} \mu,
\end{equation}
and thus the $\gamma$ for the original process $Z$ should be
\[
\gamma= U G^{-1/2} \overline{\gamma} = 2 D^{-1}
Q^{-1} \mu,
\]
just as in~(\ref{eqgamma}).

In the remainder of the proof, we go through the computations of~\cite
{HW} as applied to our setting.
Let $\overline{L}:= \frac{1}{2} \Delta+ \overline{\mu} \cdot
\nabla$,
and for $j \in[n-1 ]$, let $\overline{\mathcal{D}}_j:=
\overline{q}_{\cdot, j} \cdot\nabla$. The generator $\overline
{\mathcal
{L}}$ of the sticky Brownian motion $\overline{Z}$ can then be written as
\[
\overline{\mathcal{L}} = \overline{L} \mathbf{1}_{ \{
\overline
{\mathcal{S}} \setminus\partial\overline{\mathcal{S}} \}} + \sum
_{j=1}^{n-1} \overline{\mathcal{D}}_j
\mathbf{1}_{ \{ \overline
{F}_j \}}.
\]
Let $\overline{p} ( \overline{z} ):= \exp(
\langle\overline{\gamma}, \overline{z} \rangle)$ for
$\overline{z} \in\overline{\mathcal{S}}$. In order to show
that~(\ref{eqinvdist}) is invariant for $Z$, we must show that for
every $f \in
C_c^{\infty} ( \overline{\mathcal{S}} )$, we have
%
%
\begin{equation}
\label{eqinv} \int_{\overline{\mathcal{S}}} \overline{p} \overline{L}
f \,\mathrm
{d}\overline{z} + \sum_{j=1}^{n-1}
\frac{\sqrt{a_{j,j}}}{2} \int_{\overline{F}_j} \overline{p} \overline{
\mathcal{D}}_j f \,\mathrm{d} \overline{z}{}^j = 0,
\end{equation}
where for $j \in[n-1 ]$, $\mathrm{d}\overline{z}{}^j$ is the
surface measure on the face $\overline{F}_j$. Define $\overline{L}{}^{*}:= \frac{1}{2} \Delta- \overline{\mu} \cdot\nabla$. Using Green's
second identity and the divergence theorem, we get that for every $f
\in C_c^{\infty} ( \overline{\mathcal{S}} )$, we have
%
%
\begin{eqnarray}\label{eqGreen}
\int_{\overline{\mathcal{S}}} \overline{p} \overline
{L} f \,\mathrm
{d}\overline{z}
&=&  \int_{\overline{\mathcal{S}}} f \overline{L}{}^{*}
\overline{p} \,\mathrm{d}\overline{z}
\nonumber\\[-8pt]\\[-8pt]
&&{} + \frac{1}{2} \sum
_{j=1}^{n-1} \int_{\overline{F}_j} \biggl( f
\frac{\partial\overline
{p}}{\partial\overline{\mathfrak{n}}_{\cdot,j}} - \overline{p} \frac
{\partial f}{\partial\overline{\mathfrak{n}}_{\cdot,j}} - 2 \overline
{\mu} \cdot
\overline{\mathfrak{n}}_{\cdot,j} f \overline{p} \biggr) \,\mathrm{d}
\overline{z}{}^j,\nonumber
\end{eqnarray}
where we used $\partial/ \partial\overline{\mathfrak{n}}_{\cdot,j}
\equiv\overline{\mathfrak{n}}_{\cdot,j} \cdot\nabla$ to denote
differentiation in the inward unit normal direction on the face
$\overline{F}_j$. Now $\overline{L}{}^{*} \overline{p} = ( \frac
{1}{2} \llvert \overline{\gamma} \rrvert^2 - \overline{\mu} \cdot
\overline{\gamma} ) \overline{p} = 0$, since using~(\ref
{eqgammaformula}) we have that
\[
\tfrac{1}{2} \llvert\overline{\gamma} \rrvert^2 - \overline{
\mu} \cdot\overline{\gamma} = \tfrac{1}{2} \bigl( \overline{
\mathfrak{N}}^{-1} \overline{\gamma} \bigr)^T \overline{
\mathfrak{N}}^T \overline{\mathfrak{T}} \bigl( \overline{
\mathfrak{N}}^{-1} \overline{\gamma} \bigr),
\]
which is zero, since $\overline{\mathfrak{N}}^T \overline{\mathfrak
{T}}$ is skew symmetric due to~(\ref{eqskew}) and our condition~(\ref
{eqinvcond}). Plugging~(\ref{eqGreen}) back into~(\ref{eqinv}) and
using~(\ref{eqqdecomp}) we get that showing~(\ref{eqinv}) is
equivalent to showing that for every $f \in C_c^{\infty} (
\overline{\mathcal{S}} )$, we have
%
%
\begin{equation}
\label{eqHWbdry} \sum_{j=1}^{n-1} \int
_{\overline{F}_j} \bigl\{ f \bigl( ( \overline{\mathfrak{n}}_{\cdot,j}
- \overline{\mathfrak{t}}_{\cdot,j} ) \cdot\nabla\overline{p} - 2
\overline{\mu} \cdot\overline{\mathfrak{n}}_{\cdot, j} \overline{p}
\bigr)
+ \nabla\cdot( \overline{\mathfrak{t}}_{\cdot,j} \overline{p} f ) \bigr
\}
\,\mathrm{d}\overline{z}{}^j = 0.
\end{equation}
The relationship~(\ref{eqgammaformula}) between $\overline{\gamma}$
and $\overline{\mu}$ implies that $\overline{\gamma}{}^T (
\overline
{\mathfrak{N}} - \overline{\mathfrak{T}} ) - 2 \overline{\mu}{}^T
\overline{\mathfrak{N}} = 0$, and thus for every $j \in
[n-1
]$ we have $ ( \overline{\mathfrak{n}}_{\cdot,j} - \overline
{\mathfrak{t}}_{\cdot,j} ) \cdot\nabla\overline{p} - 2
\overline
{\mu} \cdot\overline{\mathfrak{n}}_{\cdot, j} \overline{p} = 0$.
Since $\overline{\mathfrak{t}}_{\cdot,j}$ is parallel to the face
$\overline{F}_j$, the divergence in~(\ref{eqHWbdry}) is the same as
the divergence taken in $\overline{F}_j$.
Thus, by applying the divergence theorem on each face $\overline{F}_j$,
it follows that showing~(\ref{eqHWbdry}) is equivalent to showing
that for every $f \in C_c^{\infty} ( \overline{\mathcal{S}}
)$, we have
%
%
\begin{equation}
\label{eqHWbdry2} \sum_{j=1}^{n-1} \sum
_{1 \leq k < j} \int_{\overline{F}_{j,k}} ( \overline{
\mathfrak{t}}_{\cdot,j} \cdot\overline{n}_{j,k} + \overline{
\mathfrak{t}}_{\cdot,k} \cdot\overline{n}_{k,j} ) \overline{p} f
\,\mathrm{d} \overline{\sigma}_{j,k} = 0,
\end{equation}
where $\overline{F}_{j,k} = \overline{F}_j \cap\overline{F}_k$,
$\overline{\sigma}_{j,k}$ denotes $ (n-3 )$-dimensional
surface measure on $\overline{F}_{j,k}$, and $\overline{n}_{j,k}$
denotes the unit vector that is normal to both $\overline{F}_{j,k}$ and
$\overline{\mathfrak{n}}_{\cdot,j}$, and points into the interior of
$\overline{F}_j$ from $\overline{F}_{j,k}$. In fact, $\overline
{n}_{j,k}$ must lie in the two-dimensional space spanned by $\overline
{\mathfrak{n}}_{\cdot,j}$ and $\overline{\mathfrak{n}}_{\cdot,k}$, and
can be determined uniquely,
\[
\overline{n}_{j,k} = \bigl( \overline{\mathfrak{n}}_{\cdot,k} -
\overline{\mathfrak{n}}_{\cdot,j}^T \overline{
\mathfrak{n}}_{\cdot,k} \overline{\mathfrak{n}}_{\cdot,j} \bigr) /
\bigl( 1 - \bigl( \overline{\mathfrak{n}}_{\cdot,j}^T
\overline{\mathfrak{n}}_{\cdot,k} \bigr)^2 \bigr)^{1/2}.
\]
Consequently, since $\overline{\mathfrak{N}}^T \overline{\mathfrak{T}}$
is skew symmetric, we have $\overline{\mathfrak{t}}_{\cdot,j} \cdot
\overline{n}_{j,k} + \overline{\mathfrak{t}}_{\cdot,k} \cdot
\overline
{n}_{k,j} = 0$ for all $1 \leq k < j \leq n -1$, showing that~(\ref
{eqHWbdry2}) indeed holds.
\end{pf}

\section{Convergence and general setup} \label{secconv}

This section is divided into two parts. In the first part (Section~\ref
{secconvbasic}) we prove the convergence theorem (Theorem~\ref
{mainthm}) as stated in the \hyperref[sec1]{Introduction}. Then in the second part
(Section~\ref{secgeneral}) we describe a much larger class of particle
systems that converge to appropriate sticky Brownian motions in~$\mathcal{W}$.

\subsection{Proof of the convergence theorem} \label{secconvbasic}

Given the uniqueness of a weak solution to the system of SDEs~(\ref
{mainsde}) as proved in Theorem~\ref{teouniq}, Theorem~\ref{mainthm}
is a consequence of Proposition~\ref{propconv} below. To state and
obtain the latter, we study the following decomposition. For each $i
\in[n ]$ we can write
%
%
\begin{eqnarray}
\label{Xdecomp} X_i^M ( t ) & =& X_i^M
( 0 ) + A_i^M ( t ) + \sum_{j=1}^{n-1}
C_{i,j}^{R,M} ( t ) - \sum_{j=1}^{n-1}
C_{i,j}^{L,M} ( t )
\nonumber\\[-8pt]\\[-8pt]
&&{} + \sum_{j=1}^{n-1} \Delta_{i,j}^{R,M}
( t ) - \sum_{j=1}^{n-1}
\Delta_{i,j}^{L,M} ( t ),\nonumber
\end{eqnarray}
where, for $j \in[ n - 1 ]$,
%
%
\begin{eqnarray}
A_i^M ( t ) &:=& \int_0^t
\mathbf{1}_{ \{ X_k^M
(s ) + ({1}/{\sqrt{M}}) < X_{k+1}^M ( s ), k \in
[n-1 ] \}} \,\mathrm{d} \bigl( P_i ( s ) -
Q_i ( s ) \bigr),\label{eqAi^M}
\\
C_{i,j}^{R,M} ( t ) &:=& \theta_{i,j}^R
I_{i,j}^{R,M}(t)
\nonumber\\[-8pt]\label{eqC{i,j}^{R,M}}\\[-8pt]
&:=& \theta_{i,j}^R \int_0^t
\mathbf{1}_{ \{ X_i^M ( s
)+ ({1}/{\sqrt{M}}) < X_{i+1}^M ( s ), X_j^M ( s
) + ({1}/{\sqrt{M}}) = X_{j+1}^M ( s ) \}} \,\mathrm{d} s\nonumber
\end{eqnarray}
and
%
\begin{eqnarray}\label{eqDelta{i,j}^{R,M}}
&& \Delta_{i,j}^{R,M} ( t )
\nonumber\\[-8pt]\\[-8pt]
&&\qquad := \int_0^t \mathbf{1}_{
\{X_i^M ( s ) + ({1}/{\sqrt{M}}) < X_{i+1}^M ( s
), X_j^M ( s ) + ({1}/{\sqrt{M}}) = X_{j+1}^M ( s
) \}} \,\mathrm{d}
\bigl( R_{i,j} ( s ) - \theta_{i,j}^R s \bigr),\hspace*{-18pt}\nonumber
\end{eqnarray}
and the processes $C_{i,j}^{L,M}$, $I_{i,j}^{L,M}$ and $\Delta
_{i,j}^{L,M}$ are defined similarly to $C_{i,j}^{R,M}$, $I_{i,j}^{R,M}$
and $\Delta_{i,j}^{R,M}$, respectively. For $m\in\mathbb{N}$, let
$D^m \equiv
D ([0,\infty),\mathbb{R}^m )$. We have the following
convergence result.

%
\begin{prop}\label{propconv}
Assume that Assumption~\ref{mainass} holds and that the initial
conditions $ \{ X^M (0 ), M>0 \}$ are deterministic
and converge to a limit $x\in\mathcal{W}$ as $M\rightarrow\infty$. Then
the family
%
%
\begin{equation}
\label{eqfam} \bigl\{ \bigl( X^M, A^M, I^{L,M},
I^{R,M}, \Delta^{L,M}, \Delta^{R,M} \bigr), M > 0
\bigr\}
\end{equation}
is tight in $D^{4n^2 - 2n}$. Moreover, every limit point
\[
\bigl( X^\infty, A^\infty, I^{L,\infty}, I^{R,\infty},
\Delta^{L,\infty
}, \Delta^{R,\infty} \bigr)
\]
satisfies the following for each $i \in[ n ]$:
%
%
\begin{eqnarray}
X_i^\infty( \cdot) & =& \int
_0^\cdot\mathbf{1}_{ \{ X_1^\infty( s ) <
\cdots< X_n^{\infty} ( s ) \}} \sqrt{2a}
\,\mathrm{d}W_i (s )
\nonumber\\[-8pt]\label{eqconvX}  \\[-8pt]
&&{} + \sum_{j=1}^{n-1} v_{i,j} \int
_0^\cdot\mathbf{1}_{ \{
X_j^{\infty} ( s ) = X_{j+1}^\infty( s )
\}
} \,\mathrm{d}s,\nonumber
\\
A_i^\infty( \cdot) & =& \int_0^\cdot
\mathbf{1}_{
\{
X_1^\infty( s ) < X_2^{\infty} ( s ) <
\cdots<
X_n^{\infty} ( s ) \}} \sqrt{2a} \,\mathrm{d}W_i ( s ),
\label{eqconvA}
\\
I_{i,j}^{L,\infty} ( \cdot) & =& \int_0^\cdot
\mathbf{1}_{ \{ X_j^{\infty} ( s ) = X_{j+1}^\infty( s
) \}} \,\mathrm{d}s, \qquad j \in[n-1 ] \setminus\{i-1 \},
\label{eqconvCL}
\\
I_{i,j}^{R,\infty} ( \cdot) & =& \int_0^\cdot
\mathbf{1}_{ \{ X_j^{\infty} ( s ) = X_{j+1}^\infty( s
) \}} \,\mathrm{d}s, \qquad j \in[n-1 ] \setminus\{i \},
\label{eqconvCR}
\\
I_{i,i-1}^{L,\infty} ( \cdot) & =& I_{i,i}^{R,\infty}
( \cdot) = 0,
\nonumber
\\
\Delta_{i,j}^{L,\infty} (\cdot) & =& \Delta_{i,j}^{R,\infty
}
(\cdot) = 0, \qquad j \in[n-1 ],
\nonumber
\end{eqnarray}
with a suitable $n$-dimensional standard Brownian motion $W = (
W_1,\ldots, W_n )$.
\end{prop}

\begin{pf}
\textit{Step} 1. The tightness of the family in~(\ref{eqfam}) can be
verified using the necessary and sufficient conditions of
Corollary~3.7.4 in~\cite{EK}.
First, note that for $i\in[ n ]$ and $j \in
[n-1
]$, the processes $P_i(\cdot)-Q_i(\cdot)$, as well as $ (
M^{1/4}(R_{i,j}(t)-\theta_{i,j}^R t), t\geq0 )$ and $ (
M^{1/4}(L_{i,j}(t)-\theta_{i,j}^L t), t\geq0 )$ all converge to
suitable one-dimensional Brownian motions in the limit $M\rightarrow
\infty$.
Therefore, the conditions of Corollary~3.7.4 in~\cite{EK} hold for the
corresponding families of processes indexed by $M>0$.
One can then bound the indicator functions appearing in the integrands
of the integrals in~(\ref{eqAi^M}),~(\ref{eqC{i,j}^{R,M}})
and~(\ref{eqDelta{i,j}^{R,M}}) between\vspace*{2pt} $0$ and $1$ appropriately to
show that the same conditions hold for the family
$ \{ ( A^M,I^{L,M},I^{R,M},\Delta^{L,M},\Delta^{R,M}
),
M > 0 \}$, which is thus tight in $D^{4n^2 - 3n}$.
For example, for $i \in[n ]$ and $t \geq0$ we have that
\[
\bigl\llvert A_i^M ( t ) \bigr\rrvert\leq\int
_0^t \operatorname{sgn} \bigl(
P_i ( s ) - Q_i ( s ) \bigr) \,\mathrm{d} \bigl(
P_i ( s ) - Q_i ( s ) \bigr).
\]
The expression on the right-hand side converges to $\sqrt{2a} \int_0^t
\operatorname{sgn} ( B ( s ) ) \,\mathrm
{d}B ( s )$ as
$M \to\infty$, where $B$ is a standard one-dimensional Brownian motion.
By Tanaka's formula, this is equal to $\sqrt{2a} ( \llvert
B (t
) \rrvert - L ( t ) )$, where $L (
\cdot
)$ is\vspace*{1pt} the local time process at $0$ of $B ( \cdot)$.
Consequently, the family of processes $ \{ A^M, M > 0 \}$
is tight.
Verifying tightness of the other families of processes can be done similarly.
In view of decomposition~(\ref{Xdecomp}), the first statement of the
proposition now readily follows.

\textit{Step} 2. Now fix a limit point $ ( X^\infty,
A^\infty, I^{L,\infty}, I^{R,\infty}, \Delta^{L,\infty}, \Delta
^{R,\infty} )$ and to simplify notation assume that it is the
limit of the whole family~(\ref{eqfam}) as $M \to\infty$.

We start with a few simple observations about the limit point under
consideration. Note first that, for any fixed $M>0$, the jumps of all
components of $ (X^M, A^M, I^{L,M}, I^{R,M}, \Delta^{L,M}, \Delta
^{R,M} )$ are bounded above in absolute value by $\frac
{1}{\sqrt{M}}$, so all components of the limit point must have
continuous paths.
Moreover, for every fixed $t\geq0$, the family $ \{ A^{M} ( t
), M > 0 \}$ is uniformly integrable due to the estimate
\begin{eqnarray*}
\mathbb E \bigl[A^M_i (t )^2 \bigr] & =&
\mathbb E \bigl[ \bigl[ A_i^M \bigr] ( t ) \bigr]
\\
& =& \mathbb E \biggl[\int_0^t
\mathbf{1}_{ \{ X_k^M ( s
) +
({1}/{\sqrt{M}}) < X_{k+1}^M ( s ), k\in[n-1]
\}} \,\mathrm{d} [P_i - Q_i ](s)
\biggr]
\\
& =& \mathbb E \biggl[\int_0^t
\mathbf{1}_{ \{ X_k^M ( s
) +
({1}/{\sqrt{M}}) < X_{k+1}^M ( s ), k\in[n-1]
\}} \frac{1}{\sqrt{M}} (\mathrm{d} P_i +
\mathrm{d} Q_i ) ( s ) \biggr]
\\
&\leq& 2 a t,
\end{eqnarray*}
$i \in[n ]$, where $ [\cdot]$ denotes the
quadratic variation process of a process with paths in~$D^1$. This and
the fact that $A^M$ is a martingale for any fixed $M>0$ show that
$A^{\infty}$ is a martingale with respect to its own filtration; see,
for example,~\cite{js}, Proposition~IX.1.12.

Next, we observe that, as limits of nondecreasing processes,
$I_{i,j}^{L,\infty}$ and $I_{i,j}^{R,\infty}$ must be nondecreasing
processes themselves for every $i\in[n ]$, $j \in
[n-1 ]$, and consequently they are also of finite variation.
Furthermore, for all $i\in[n ]$, $j\in[n-1
]$, the
quadratic variation processes of the martingales $\Delta_{i,j}^{L,M}$
and $\Delta_{i,j}^{R,M}$ satisfy
\[
\forall t\geq0: \qquad\lim_{M\rightarrow\infty} \mathbb E \bigl[ \bigl
[\Delta
_{i,j}^{L,M} \bigr](t) \bigr]=0\quad\mbox{and}\quad \lim
_{M\rightarrow\infty} \mathbb E \bigl[ \bigl[\Delta_{i,j}^{R,M}
\bigr](t) \bigr]=0.
\]
Therefore, the distributional limits
\[
\Delta_{i,j}^{L,\infty}\equiv\lim_{M\rightarrow\infty} \Delta
_{i,j}^{L,M},\qquad\Delta_{i,j}^{R,\infty}\equiv
\lim_{M\rightarrow\infty} \Delta_{i,j}^{R,M}
\]
in $D^1$ exist and are identically equal to zero.

Finally, the two observations of the previous paragraph, together with
the decomposition~(\ref{Xdecomp}), show that for $i,i'\in[n]$ the
quadratic covariation processes $ \langle X_i^{\infty},
X_{i'}^\infty\rangle$ and $ \langle A_i^{\infty
},A_{i'}^\infty
\rangle$ are in fact equal. In particular, we have that $
\langle X_i^\infty\rangle= \langle A_i^\infty
\rangle
$ for $i \in[ n ]$.

\textit{Step} 3. In order to show~(\ref{eqconvA}), we study
the quadratic covariation processes $ \langle X_i^{\infty
},X_{i'}^{\infty} \rangle= \langle A_i^{\infty
},A_{i'}^\infty
\rangle$, $i,i'\in[n]$.
We first claim that $ \langle A_i^{\infty},A_{i'}^\infty
\rangle=0$ whenever $i\neq i'$. To this end, it suffices to show that
for any such pair of indices $A_i^{\infty}(\cdot)A_{i'}^\infty(\cdot)$
is a martingale with respect to its own filtration. The latter is the
limit in $D^1$ of the family of martingales $ \{ A_i^M(\cdot
)A_{i'}^M(\cdot), M>0 \}$ by definition, so it is enough to prove
that, for any fixed $t\geq0$, the random variables $ \{
A_i^M(t)A_{i'}^M(t), M>0 \}$ are uniformly integrable. The latter
is a consequence of the following chain of estimates:
%
\begin{eqnarray*}
&&\mathbb E \bigl[A_i^M(t)^2
A_{i'}^M(t)^2 \bigr]
\\
&&\qquad  =\mathbb E \biggl[\int_0^t
A_i^M(s)^2 \,\mathrm{d}A^M_{i'}(s)^2
\biggr]+\mathbb E \biggl[\int_0^t
A^M_{i'}(s)^2 \,\mathrm{d}A^M_i(s)^2
\biggr]
\\
&&\qquad  =\mathbb E \biggl[\int_0^t
A_i^M(s)^2 \,\mathrm{d} \bigl[A^M_{i'}
\bigr](s) \biggr]+\mathbb E \biggl[\int_0^t
A_{i'}^M(s)^2 \,\mathrm{d} \bigl[A^M_i
\bigr](s) \biggr]
\\
&&\qquad \leq\mathbb E \biggl[\int_0^t
A_i^M(s)^2 \,\mathrm{d} [P_{i'}-Q_{i'}
](s) \biggr]+\mathbb E \biggl[\int_0^t
A_{i'}^M(s)^2 \,\mathrm{d} [P_i-Q_i
](s) \biggr]
\\
&&\qquad  =\mathbb E \biggl[\int_0^t
A_i^M(s)^2 \frac{1}{\sqrt{M}}(\mathrm
{d}P_{i'}+\mathrm{d}Q_{i'}) (s) \biggr]
\\
&&\qquad\quad{} +\mathbb E \biggl[\int_0^t
A_{i'}^M(s)^2 \frac{1}{\sqrt{M}}(\mathrm
{d}P_i+\mathrm{d}Q_i) (s) \biggr]
\\
&&\qquad  =\mathbb E \biggl[\int_0^t
A_i^M(s)^2 2a \,\mathrm{d}s \biggr]+\mathbb E
\biggl[\int_0^t A_{i'}^M(s)^2
2a \,\mathrm{d}s \biggr]
\\
&&\qquad \leq2a \int_0^t \bigl( \mathbb E \bigl[
[P_i-Q_i ](s) \bigr]+\mathbb E \bigl[
[P_{i'}-Q_{i'} ](s) \bigr] \bigr) \,\mathrm{d}s
\\
&&\qquad  = 2a \int_0^t 4as \,\mathrm{d}s=4a^2t^2.
\end{eqnarray*}

We next aim to evaluate $ \langle A_i^\infty\rangle$,
$i\in
[n ]$. To this end, we first show that for any fixed $i\in
[n]$ and $t\geq0$, the random variables $ \{ A^M_i(t)^2, M>0
\}$ are uniformly integrable. By It\^{o}'s lemma for (not necessarily
continuous) semimartingales we have that for any fixed $t \geq0$,
\begin{eqnarray*}
A_i^M ( t )^4 & =& \int
_0^t 4 A_i^M ( s-
)^3 \,\mathrm{d} A_i^M ( s ) +
\frac{1}{2} \int_0^t 12
A_i^M ( s- )^2 \,\mathrm{d} \bigl[
A_i^M \bigr] ( s )
\\
&&{} + \sum_{s \leq t} \bigl\{ \bigl( A_i^M
( s )^4 - A_i^M ( s- )^4
\bigr) - 4 A_i^M ( s- )^3 \bigl(
A_i^M ( s ) - A_i^M ( s- )
\bigr)
\\
&&\hspace*{125pt}{}  - 6 A_i^M ( s- )^2 \bigl(
A_i^M ( s ) - A_i^M ( s- )
\bigr)^2 \bigr\}.
\end{eqnarray*}
Taking expectations on both sides and dropping the third line of the
previous display, we arrive at the following estimate:
\begin{eqnarray*}
\mathbb E \bigl[A_i^M(t)^4 \bigr] &\leq&
\mathbb E \biggl[\int_0^t 4
A_i^M(s-)^3 \,\mathrm{d}A_i^M(s)
\biggr]
\\
&&{}  + \mathbb E \biggl[\int_0^t 6
A_i^M(s-)^2 \,\mathrm{d}
\bigl[A_i^M \bigr] ( s ) \biggr]
\\
&&{} + \mathbb E \biggl[ \sum_{s \leq t} \bigl\{ \bigl(
A_i^M ( s )^4 - A_i^M
( s- )^4 \bigr)
\\
&&\hspace*{42pt}{} - 4 A_i^M ( s-
)^3 \bigl( A_i^M ( s ) -
A_i^M ( s- ) \bigr) \bigr\} \biggr].
\end{eqnarray*}
The first term on the right-hand side is equal to zero, since $A_i^M$
is a martingale starting at zero. We can upper bound the second term on
the right-hand side using the inequality $\mathrm{d} [A_i^M
]
( s ) \leq\frac{1}{\sqrt{M}}\,\mathrm{d} ( P_i + Q_i
) ( s )$. Finally, we can upper bound the third
expectation on the right-hand side using the identity $x^4 - y^4 - 4y^3
( x - y ) = ( x - y )^2 ( x^2 + 2xy + 3y^2
)$, the fact that the jumps of the process $A_i^M$ are of size
$\frac{1}{\sqrt{M}}$, and the fact that $A_i^M ( s ) - A_i^M
( s- ) \leq\mathrm{d} ( P_i + Q_i ) ( s
)$. We then arrive at the estimate
\[
\mathbb E \bigl[A_i^M(t)^4 \bigr] \leq
\mathbb E \biggl[\int_0^t \bigl( 3
A_i^M(s)^2 + 11 A_i^M
( s- )^2 \bigr) \frac{1}{\sqrt{M}} \,\mathrm{d}(P_i+Q_i)
(s) \biggr],
\]
and the right-hand side can now be bounded above by a constant
depending only on $a$ and $t$ by arguing as above when estimating
$\mathbb E
[A_i^M(t)^2 A_{i'}^M(t)^2 ]$. Thus indeed the random
variables $ \{ A^M_i(t)^2, M>0 \}$ are uniformly integrable.
Putting this together with the fact that the functional
\[
(\omega_1,\omega_2,\ldots,\omega_n )
\mapsto\int_{t_1}^{t_2} \mathbf{1}_{ \{ \omega_1 ( s ) <
\cdots<
\omega_n ( s ) \}}
\,\mathrm{d}s
\]
on $D^n$ is lower semicontinuous and using the Portmanteau theorem, we
have for each $i \in[n ]$ that
\begin{eqnarray*}
\hspace*{-6pt}&& \biggl[G \bigl(A^\infty\bigr) \biggl( A_i^{\infty}
(t_2 )^2 - A_i^{\infty}
(t_1 )^2 - 2a \int_{t_1}^{t_2}
\mathbf{1}_{\{ X_1^\infty(s)<\cdots
<X_n^{\infty
}(s)\}} \,\mathrm{d}s \biggr) \biggr]
\\
\hspace*{-6pt}&&\!\qquad  =\lim_{\epsilon\downarrow0} \mathbb E \biggl[G \bigl(A^\infty
\bigr) \biggl( A_i^{\infty} (t_2 )^2 -
A_i^{\infty} (t_1 )^2
\\
\hspace*{-6pt}&&\!\hspace*{67pt}\quad\qquad{} - 2a \int_{t_1}^{t_2}
\mathbf{1}_{\{ X_k^\infty(s)+\epsilon
<X_{k+1}^\infty(s), k\in[n-1]\}} \,\mathrm{d}s \biggr) \biggr]
\\
\hspace*{-6pt}&&\!\qquad \geq\limsup_{M\rightarrow\infty} \mathbb E \biggl[G \bigl(A^M
\bigr) \biggl( A_i^M (t_2 )^2 -
A_i^M (t_1 )^2
\\
\hspace*{-6pt}&&\!\hspace*{83pt}\quad\qquad{} - 2a\int_{t_1}^{t_2}
\mathbf{1}_{\{ X^M_k(s)+({1}/{\sqrt
{M}})<X^M_{k+1}(s),k\in[n-1]\}}\,\mathrm{d}s \biggr) \biggr]
\\
\hspace*{-6pt}&&\!\qquad  =\limsup_{M\rightarrow\infty} \mathbb E \biggl[G \bigl(A^M
\bigr) \biggl( A_i^M (t_2 )^2 -
A_i^M (t_1 )^2
\\
\hspace*{-6pt}&&\!\hspace*{84pt}\quad\qquad{} - \int_{t_1}^{t_2}
\mathbf{1}_{\{ X^M_k(s)+({1}/{\sqrt
{M}})<X^M_{k+1}(s),k\in[n-1]\}}\,\mathrm{d}[P_i-Q_i](s) \biggr)
\biggr]
\\
\hspace*{-6pt}&&\!\qquad  = \limsup_{M\rightarrow\infty} \mathbb E \bigl[G \bigl(A^M
\bigr) \bigl( A_i^M (t_2 )^2 -
A_i^M (t_1 )^2
- \bigl[A_i^M \bigr](t_2)
+ \bigl[A_i^M \bigr](t_1) \bigr) \bigr]=0
\end{eqnarray*}
for any nonnegative continuous bounded functional $G$ on $D^n$
measurable with respect to the $\sigma$-algebra generated by the
coordinate mappings on $D^n ( [0,t_1 ] )$.
Therefore, recalling from the end of step~2 of the proof that $
\langle X_i^\infty\rangle= \langle A_i^\infty
\rangle
$, we conclude that
%
%
\begin{equation}
\label{eqport1} \forall0 \leq t_1 < t_2\dvtx \qquad\bigl
\langle X_i^{\infty} \bigr\rangle( t_2 ) - \bigl
\langle X_i^{\infty} \bigr\rangle( t_1 ) \geq2a
\int_{t_1}^{t_2} \mathbf{1}_{ \{ X_1^\infty( s
)
< \cdots< X_n^{\infty} ( s ) \}}
\,\mathrm{d}s\hspace*{-25pt}
\end{equation}
holds with probability one. On the other hand,
\begin{eqnarray*}
&&\mathbb E \bigl[G \bigl(A^\infty\bigr) \bigl( A_i^{\infty}
(t_2 )^2 - A_i^{\infty}
(t_1 )^2 - 2a (t_2-t_1) \bigr)
\bigr]
\\
&&\qquad=\lim_{M\rightarrow\infty} \mathbb E \bigl[G \bigl(A^M
\bigr) \bigl( A_i^M (t_2 )^2 -
A_i^M (t_1 )^2 - 2a
(t_2-t_1) \bigr) \bigr]
\\
&&\qquad=\lim_{M\rightarrow\infty} \mathbb E \bigl[G \bigl(A^M
\bigr) \bigl( A_i^M (t_2 )^2 -
A_i^M (t_1 )^2
- [P_i-Q_i](t_2)
+ [P_i-Q_i](t_1) \bigr) \bigr]
\\
&&\qquad \leq 0
\end{eqnarray*}
for any functional $G$ on $D^n$ as above. Hence
%
%
\begin{equation}
\label{eqport2} \forall0\leq t_1 < t_2\dvtx \qquad\bigl
\langle X_i^{\infty} \bigr\rangle( t_2 ) - \bigl
\langle X_i^{\infty} \bigr\rangle( t_1 ) \leq2a (
t_2 - t_1 )
\end{equation}
must hold with probability 1.

In view of (\ref{eqport2}), we see that in order to improve (\ref
{eqport1}) to an equality, it suffices to show that the measure
$\mathrm
{d} \langle X^\infty_i \rangle=\mathrm{d} \langle
A^\infty
_i \rangle$ assigns zero mass to the sets $\{t\geq0\dvtx  X^\infty
_j(t)=X^\infty_{j+1}(t)\}$, $j\in[n-1]$, with probability one. To this
end, we first recall that for every $i \in[n ]$ the square
integrable martingale $A^\infty_i$ is the limit in $D^1$ of the square
integrable martingales $ \{ A^M_i, M>0 \}$, and the random
variables $ \{ A^M_i(t)^2, M>0 \}$ are uniformly integrable
for any fixed $t\geq0$. Therefore $ \langle A^\infty_i
\rangle
$ is the limit in $D^1$ of $ \{ [A^M_i ], M>0 \}$,
and so by the Portmanteau theorem,
\begin{eqnarray*}
&&\mathbb E \bigl[ \bigl( \bigl\langle A^\infty_i \bigr
\rangle(t)- \bigl\langle A^\infty_{i'} \bigr\rangle(t)
\bigr)^2 \bigr]
\\
&&\qquad \leq\liminf_{M\rightarrow\infty} \mathbb E \bigl[
\bigl( \bigl[A^M_i \bigr](t)- \bigl[A^M_{i'}
\bigr](t) \bigr)^2 \bigr]
\\
&&\qquad\leq\liminf_{M\rightarrow\infty} \mathbb E \biggl[ \biggl(\int
_0^t \mathbf{1}_{\{ X^M_k(s)+({1}/{\sqrt{M}})<X^M_{k+1}(s),k\in[n-1]\}
}
\\
&&\hspace*{94pt}{}\times \frac{1}{\sqrt{M}}
(\mathrm{d}P_i+\mathrm
{d}Q_i-\mathrm{d}P_{i'}-\mathrm{d}Q_{i'} ) (s)
\biggr)^2 \biggr]
\\
&&\qquad= \liminf_{M\rightarrow\infty} \mathbb E \biggl[\int
_0^t \mathbf{1}_{\{ X^M_k(s)+({1}/{\sqrt{M}})<X^M_{k+1}(s),k\in
[n-1]\}}
\\
&&\hspace*{88pt}{}\times
\frac
{1}{M}\,
\mathrm{d} [P_i+Q_i-P_{i'}-Q_{i'}
](s) \biggr]
\\
&&\qquad\leq\liminf_{M\rightarrow\infty} \frac{1}{M^{3/2}} \mathbb E
\bigl[P_i(t)+Q_i(t)+P_{i'}(t)+Q_{i'}(t)
\bigr] = 0
\end{eqnarray*}
for any fixed $i,i'\in[n]$ and $t\geq0$ with probability one. In view
of the path continuity of the processes $ \langle X^\infty
_i
\rangle$, $i\in[n]$, this implies
\[
\bigl\langle X^\infty_1 \bigr\rangle= \bigl\langle
X^\infty_2 \bigr\rangle= \cdots= \bigl\langle
X^\infty_n \bigr\rangle
\]
with probability one.
To conclude the argument, we use the occupation time formula for
continuous semimartingales (see, e.g.,~\cite{RY}, Theorem~VI.1.6),
which states that if $Y$ is a continuous semimartingale, and $\phi$ is
a positive Borel function, then a.s. for every $t \geq0$ we have
\[
\int_0^t \phi\bigl( Y ( s ) \bigr) \,\mathrm
{d} \langle Y \rangle( s ) = \int_{-\infty}^{\infty
} \phi
( a ) L^a ( t ) \,\mathrm{d}a,
\]
where $L^a ( \cdot)$ is the local time process at $a$ of $Y
( \cdot)$. In particular, the choice of $\phi( a
) = \mathbf{1}_{ \{ a = 0 \}}$ gives that
\[
\int_0^t \mathbf{1}_{ \{ Y ( s ) = 0 \}}
\,\mathrm{d} \langle Y \rangle( s ) = 0,
\]
and now choosing $Y ( \cdot) = X_{j+1}^{\infty} (
\cdot) - X_j^\infty( \cdot)$ this implies that
the measure
\[
\mathrm{d} \bigl\langle X_{j+1}^{\infty} - X_j^{\infty}
\bigr\rangle= \mathrm{d} \bigl\langle X_{j+1}^{\infty} \bigr
\rangle+\mathrm{d} \bigl\langle X_j^{\infty} \bigr\rangle=2
\,\mathrm{d} \bigl\langle X_j^{\infty} \bigr\rangle=2 \,\mathrm{d}
\bigl\langle X_i^{\infty} \bigr\rangle
\]
assigns zero mass to the set $\{t\geq0\dvtx  X^\infty_j(t)=X^\infty
_{j+1}(t)\}$ with probability one.
Hence, equality must hold in~(\ref{eqport1}).
The representation~(\ref{eqconvA}) with a suitable standard Brownian
motion $W=(W_1,W_2,\ldots,W_n)$ now readily follows from the Martingale
Representation theorem in the form of Theorem 4.2 in Chapter~3 of~\cite{KS}.

\textit{Step} 4. We now turn to the proof of~(\ref
{eqconvX}),~(\ref{eqconvCL}) and~(\ref{eqconvCR}). To this end,
recalling the ghost particles $X_0^M ( \cdot) \equiv-
\infty$ and $X_{n+1}^M ( \cdot) \equiv\infty$ introduced
for notational convenience, for any $M > 0$, $i \in\{ 0,1,\ldots,
n \}$ and $j \in[n-1 ]$, define
\begin{eqnarray*}
I_j^{M,1} ( \cdot) &:=& \int_0^\cdot
\mathbf{1}_{
\{
X_j^M ( s ) + ({1}/{\sqrt{M}}) = X_{j+1}^M ( s
) \}} \,\mathrm{d} s,
\\
I_{i,j}^{M,2} ( \cdot) &:=& \int_0^\cdot
\mathbf{1}_{
\{ X_{i}^M ( s ) + ({1}/{\sqrt{M}}) = X_{i+1}^M
( s
),
X_j^M ( s ) + ({1}/{\sqrt{M}}) = X_{j+1}^M ( s
) \}} \,\mathrm{d} s.
\end{eqnarray*}
Then for any $i \in[n ]$, $j \in[n-1 ]$, we have
the decompositions
\begin{eqnarray*}
I_{i,j}^{L,M} ( \cdot) &:=& I_j^{M,1} (
\cdot) - I_{i-1,j}^{M,2} ( \cdot),
\\
I_{i,j}^{R,M} ( \cdot) &:=& I_j^{M,1} (
\cdot) - I_{i,j}^{M,2} ( \cdot).
\end{eqnarray*}
It is now easy to check that for each $i \in\{ 0,1,\ldots,
n
\}$, $j\in[n-1 ]$, the families of processes $ \{
I^{M,1}_j, M>0 \}$ and $ \{ I^{M,2}_{i,j}, M>0 \}$
satisfy the\vspace*{2pt} tightness criterion of Corollary~3.7.4 in~\cite{EK}. So,\vspace*{1pt}
after passing to a subsequence if necessary, we obtain the existence of
suitable limits $I^{\infty,1}_j$ and $I^{\infty,2}_{i,j}$,
respectively;\vspace*{2pt} for notational convenience we assume that the full
families of processes converge jointly to the respective limit points.

The limiting processes inherit many properties of the prelimit processes.
First, clearly the limits are nondecreasing processes and inherit the
property that for every $i \in\{0,1,\ldots, n \}$,
%
%
\begin{equation}
\label{bdry1} \forall0\leq t_1<t_2\dvtx \qquad
I^{\infty,2}_{i,j}(t_2)-I^{\infty,2}_{i,j}(t_1)
\leq I^{\infty,1}_j(t_2)-I^{\infty,1}_j(t_1).
\end{equation}
Second, the prelimit processes satisfy
%
%
\begin{eqnarray}\label{eqbdry2prelim}
\int_0^\infty\mathbf{1}_{ \{X_j^M (t ) +
({1}/{\sqrt{M}}) < X_{j+1}^M ( t ) \}} \,\mathrm
{d}I_j^{M,1} ( t ) &=& 0,
\nonumber\\[-8pt]\\[-8pt]
\int_0^\infty( \mathbf{1}_{ \{X_{i}^M (t
) +
({1}/{\sqrt{M}}) < X_{i+1}^M ( t ) \}} +
\mathbf{1}_{ \{X_j^M (t ) + ({1}/{\sqrt{M}}) < X_{j+1}^M
( t ) \}} ) \,\mathrm{d} I_{i,j}^{M,2} ( t ) &=&0,\nonumber
\end{eqnarray}
and from these we have that the limiting processes satisfy
%
%
\begin{eqnarray}
\int_0^\infty\mathbf{1}_{ \{X^\infty_j(t)<X^\infty
_{j+1}(t) \}}
\,\mathrm{d}I^{\infty,1}_j(t)&=&0, \label{bdry2}
\\
\int_0^\infty( \mathbf{1}_{ \{X^\infty
_{i}(t)<X^\infty
_{i+1}(t) \}} +
\mathbf{1}_{ \{X^\infty_j(t)<X^\infty
_{j+1}(t) \}} ) \,\mathrm{d}I^{\infty,2}_{i,j}(t)&=&0.
\label{bdry3}
\end{eqnarray}
These properties can be shown by arguing as in the second half of the
proof of Theorem~4.1 in~\cite{Wi2} (see also the proof of Proposition~9
in~\cite{karatzas2012systems}); we provide a sketch on how to
obtain~(\ref{bdry2}) from~(\ref{eqbdry2prelim}), and~(\ref{bdry3})
follows similarly.
We first use the Skorokhod representation theorem~\cite{EK}, Theorem~3.1.8, and the fact that the limiting processes $ (
X^\infty, I^{\infty,1}, I^{\infty,2}, I^{L,\infty}, I^{R,\infty}
)$ are a.s. continuous to replace the sequence of processes $ \{
( X^M, I^{M,1}, I^{M,2}, I^{L,M}, I^{R,M} ), M > 0
\}$
by one that has the same distribution and which a.s. converges
uniformly on compact time intervals.
Let $ \{ f_m \}_{m \geq1}$ be a sequence of continuous
functions such that for every $m$, $f_m\dvtx  \mathbb{R}\to[ 0, 1
]$,
$f_m ( x ) = 0$ for $x \leq1/m$ and $f_m ( x
) =
1$ for $x \geq2/m$. By passing to the $m \to\infty$ limit, in order
to show~(\ref{bdry2}) it suffices to show that for each $t \geq0$, $j
\in[n-1 ]$, and $m \geq1$, a.s.
%
%
\begin{equation}
\label{eqbdry2approx} \int_0^t f_m
\bigl( X_{j+1}^\infty( s ) - X_j^\infty(
s ) \bigr) \,\mathrm{d}I_j^{\infty,1} ( s ) = 0.
\end{equation}
To do this, fix $j \in[ n -1 ]$, $m \geq1$ and $t \geq0$.
For $M > m^2$,~(\ref{eqbdry2prelim}) implies that a.s.
\[
\int_0^t f_m \bigl(
X_{j+1}^M ( s ) - X_j^M ( s )
\bigr) \,\mathrm{d}I_j^{M,1} ( s ) = 0,
\]
and thus to show~(\ref{eqbdry2approx}) it suffices to show that
a.s. %
%
%
\begin{eqnarray}
\label{eqbdry2conv} && \int_0^t f_m
\bigl( X_{j+1}^M ( s ) - X_j^M ( s
) \bigr) \,\mathrm{d}I_j^{M,1} ( s )
\nonumber\\[-8pt]\\[-8pt]
&&\qquad \to\int_0^t f_m \bigl(
X_{j+1}^\infty( s ) - X_j^\infty( s )
\bigr) \,\mathrm{d}I_j^{\infty,1} ( s )\nonumber
\end{eqnarray}
as $M \to\infty$. The almost sure convergence assumed above implies
that a.s. as \mbox{$M \to\infty$}, $X_{j+1}^M ( \cdot) - X_j^M
( \cdot) \to X_{j+1}^\infty( \cdot) -
X_j^\infty( \cdot)$ uniformly on compacts, and since $f_m$
is uniformly continuous, we have that a.s. as $M \to\infty$, $f_m
( X_{j+1}^M ( \cdot) - X_j^M ( \cdot)
) \to f_m ( X_{j+1}^\infty( \cdot) -
X_j^\infty
( \cdot) )$ uniformly\vspace*{1pt} on compacts. We also have that
a.s. as $M \to\infty$, $I_j^{M,1} ( \cdot) \to
I_j^{\infty,1} ( \cdot)$ uniformly on\vspace*{2pt} compacts, and the remaining
details of showing~(\ref{eqbdry2conv}) are as in the end of the proof
of Theorem~4.1 in~\cite{Wi2}.

Next, we define the time change
\[
\sigma(t)=\inf\biggl\{s\geq0\dvtx  \int_0^s
\mathbf{1}_{\{X^\infty
_1(r)<X^\infty_2(r)<\cdots<X^\infty_n(r)\}} \,\mathrm{d}r=t \biggr\},\qquad t\geq0
\]
and then let $\widehat{X}{}^\infty(\cdot)=X^{\infty}
(
\sigma(\cdot) )$, $\widehat{I}{}^{\infty,1}
( \cdot
) = I^{\infty,1} ( \sigma( \cdot)
)$ and
$\widehat{I}{}^{\infty,2} ( \cdot) = I^{\infty,2} (
\sigma( \cdot) )$. Using L\'{e}vy's characterization
of Brownian motion, we conclude that the components of $\widehat
{X}{}^\infty$ admit the decomposition
\begin{eqnarray*}
\widehat{X}{}^\infty_i(\cdot) & =& \widehat{X}{}^\infty_i(0)+
\sqrt{2a} \widehat{W}_i(\cdot) + \sum_{j=1}^{n-1}
v_{i,j} \widehat{I}{}^{\infty,1}_j( \cdot)
\\
&&{} + \mathop{\sum_{j=1}}_{j \neq i-1}^{n-1}
\theta^L_{i,j} \widehat{I}{}^{\infty,2}_{i-1,j} (
\cdot) - \mathop{\sum_{j=1}}_{j
\neq
i}^{n-1}
\theta^R_{i,j} \widehat{I}{}^{\infty,2}_{i,j} (
\cdot)
\end{eqnarray*}
with $\widehat{W}=(\widehat{W}_1,\widehat{W}_2,\ldots,\widehat{W}_n)$
being a suitable standard Brownian motion.
As we shall show shortly, for every $i \in[n ]$ we have
%
%
\begin{equation}
\label{eqdouble} \widehat{I}_{i,j}^{\infty,2} ( \cdot) \equiv0,
\qquad j \in[n-1 ] \setminus\{ i \},
\end{equation}
and thus the decomposition simplifies to
\[
\widehat{X}{}^\infty_i(\cdot) = \widehat{X}{}^\infty_i(0)+
\sqrt{2a} \widehat{W}_i(\cdot) + \sum_{j=1}^{n-1}
v_{i,j} \widehat{I}{}^{\infty,1}_j(\cdot).
\]
We can then identify the process of spacings
\[
\bigl(\widehat{X}{}^\infty_2(\cdot)-\widehat{X}_1^\infty(
\cdot), \widehat{X}{}^\infty_3(\cdot)-\widehat{X}_2^\infty(
\cdot), \ldots, \widehat{X}{}^\infty_n(\cdot)-
\widehat{X}_{n-1}^\infty(\cdot) \bigr)
\]
as a reflected Brownian motion in the orthant $(\mathbb{R}_+)^{n-1}$ with
reflection matrix $Q$ (recall from Section~\ref{secintroconv} that
$q_{j,j'} = v_{j+1,j'} - v_{j,j'}$ for $j,j' \in[n-1 ]$),
and the processes $\widehat{I}{}^{\infty,1}_j(\cdot)$, $j\in[n-1]$, with
its boundary local times.
At this point one can argue as in step~2 in the proof of Theorem~\ref
{teouniq} to obtain the representations~(\ref{eqconvX}),~(\ref
{eqconvCL}) and~(\ref{eqconvCR}).

Thus what is left is to show~(\ref{eqdouble}). For $j \in
[n-1 ]$ let $\widehat{Z}_j ( \cdot) = \widehat
{X}_{j+1}^\infty( \cdot) - \widehat{X}_j^\infty(
\cdot)$, thus $\widehat{Z} ( \cdot) = (
\widehat
{Z}_1 ( \cdot),\ldots, \widehat{Z}_{n-1} ( \cdot
) )$ is the process of spacings. Due to~(\ref{bdry3}),
showing~(\ref{eqdouble}) reduces to showing that
%
%
\begin{equation}
\label{eqbdryend1} \int_0^\infty\mathbf{1}_{ \{\widehat{Z}_{i}(t) =
\widehat
{Z}_{j}(t) = 0 \}}
\,\mathrm{d}\widehat{I}{}^{\infty,2}_{i,j}(t)=0.
\end{equation}
This can be done by generalizing the proof of Theorem~1 in~\cite{RW},
along the lines of~\cite{bhardwaj2009diffusion}, Theorem~7.7, and~\cite{karatzas2012systems}, Lemma~1, and, in particular, it uses the
Lyapunov functions constructed in the proof of Lemma~4 in~\cite{RW}.
Here we provide a sketch of the proof, and refer to~\cite{RW,bhardwaj2009diffusion} and~\cite{karatzas2012systems} for details.
This is the only point in our proof where we use Assumption~\ref{mainass}(b).

First we introduce some notation to simplify the representation of
$\widehat{Z}$. For $i \in[ n - 1 ]$, let $\widehat{B}_i
( \cdot
):= \sqrt{2a} ( \widehat{W}_{i+1} ( \cdot
) - \widehat{W}_i ( \cdot) )$; then $\widehat{B}:= ( \widehat{B}_1,\ldots, \widehat{B}_{n-1} )$ is a Brownian motion with mean
zero and
diffusion matrix $A = ( a_{i,j} )_{i,j = 1}^{n-1}$ given by
\[
a_{i,j}:= \cases{ 4a, &\quad if $i = j$,
\vspace*{5pt}\cr
-2a, &\quad if $\llvert i
- j \rrvert= 1$,
\vspace*{5pt}\cr
0, &\quad otherwise.}
\]
In the following we think of $\widehat{I}{}^{\infty,2}$ as an $ (
\mathbb{R}_+
)^{ (n-1 )^2}$-valued process whose components are
indexed by ordered pairs $ ( i,j )$, $i,j \in
[n-1
]$, and the component indexed by $ ( i,j )$ is $\widehat
{I}{}^{\infty,2}_{i,j} ( \cdot)$.
Recalling\vspace*{2pt} the definition of the matrix $Q^{ (2 )}$ from
Section~\ref{secintroconv}, we can write $\widehat{Z}$ as
\[
\widehat{Z} ( \cdot) = \widehat{Z} ( 0 ) + \widehat{B} ( \cdot) + Q
\widehat{I}{}^{\infty,1} ( \cdot) + Q^{
(2 )} \widehat{I}{}^{\infty,2}
( \cdot).
\]
Then by It\^o's formula, for any function $f$ that is twice
continuously differentiable in some domain containing $ ( \mathbb{R}_+
)^{n-1}$ we have that a.s. for all $t \geq0$,
\begin{eqnarray*}
&&f \bigl( \widehat{Z} ( t ) \bigr) - f \bigl( \widehat{Z} ( 0 ) \bigr)
\\
&&\qquad= \int_0^t \nabla f \bigl(
\widehat{Z} ( s ) \bigr) \,\mathrm{d}\widehat{B} ( s ) + \sum
_{j=1}^{n-1} \int_0^t
q_{\cdot,j} \cdot\nabla f \bigl( \widehat{Z} ( s ) \bigr) \,\mathrm{d}
\widehat{I}_j^{\infty,1} ( s )
\\
&&\quad\qquad{} + \sum_{k,\ell= 1}^{n-1} \int
_0^t q^{ (2 )}_{\cdot,
(k,\ell)} \cdot
\nabla f \bigl( \widehat{Z} ( s ) \bigr) \,\mathrm{d} \widehat{I}_{k,\ell
}^{\infty,2}
( s ) + \int_0^t Lf \bigl( \widehat{Z} ( s )
\bigr) \,\mathrm{d}s,
\end{eqnarray*}
where recall that $q_{\cdot,j}$ is the $j$th column of $Q$, $q^{
(2 )}_{\cdot, (k,\ell)}$ is the column of $Q^{
(2 )}$ corresponding to index $ (k,\ell)$ and
\[
L = \frac{1}{2} \sum_{i,j=1}^{n-1}
a_{i,j} \frac{\partial
^2}{\partial
x_i\, \partial x_j}.
\]
We apply It\^o's formula to an appropriately defined family of
functions, just as in~\cite{RW}.
Let $\gamma= \gamma( [n-1 ] ) \in(
\mathbb{R}_+
)^{n-1}$ be the vector\vspace*{1pt} guaranteed by Assumption~\ref{mainass}(b)
for $J = [n-1 ]$.
Let $\delta:= Q^T \gamma$; by assumption $\delta\in[1,\infty)^{n-1}$.
Define $\alpha= A \gamma$.
For each $x \in( \mathbb{R}_+ )^{n-1}$ and $r \in
( 0, 1
)$, let $d^2 ( x, r ):= ( x + r \alpha
)^T
A^{-1} ( x + r \alpha)$. Then, for each $\varepsilon\in
( 0,
1 )$, define
%
%
\begin{equation}
\label{eqphi}
\phi_{\varepsilon} ( x ):=
\cases{ \displaystyle\frac{1}{2- ( n - 1
)} \int
_{\varepsilon}^1 r^{ (n -1 )
- 2} \bigl( d^2
( x, r ) \bigr)^{({2 - ( n - 1
)})/{2}} \,\mathrm{d}r,
\cr
\hspace*{105pt}\quad\mbox{if $n - 1 \geq3$,}
\cr
\displaystyle\frac{1}{2} \int_{\varepsilon}^1 \ln\bigl(
d^2 ( x,r ) \bigr) \,\mathrm{d}r,\qquad\mbox{if $n - 1 = 2$}.}
\end{equation}
For each $\varepsilon\in( 0, 1 )$, $\phi_\varepsilon$
is twice
continuously differentiable in some domain containing $ ( \mathbb{R}_+
)^{n-1}$, and on each compact subset of $ ( \mathbb{R}_+
)^{n-1}$ it is bounded, uniformly in $\varepsilon$. Moreover, we have
that $L
\phi_\varepsilon= 0$ in some domain containing $ ( \mathbb{R}_+
)^{n-1}$,
due to the fact that the integrands in~(\ref{eqphi}) are $L$-harmonic
functions of $x \in\mathbb{R}^{n-1} \setminus\{ - r \alpha
\}$.
Now, with $\llVert \cdot\rrVert$ denoting the Euclidean norm in
$\mathbb{R}
^{n-1}$, define for each $m \in\mathbb{N}$ the stopping time
\[
\tau_m:= \inf\bigl\{ t \geq0\dvtx  \bigl\llVert\widehat{Z} ( t )
\bigr\rrVert\geq m\mbox{ or }\widehat{I}_j^{\infty,1} (
t ) \geq m\mbox{ for some } j \bigr\} \wedge m.
\]
Applying It\^o's formula to the function $\phi_\varepsilon$ and the stopping
time $\tau_m$, we get that a.s. %
%
%
\begin{eqnarray}
\label{eqitophi} \phi_{\varepsilon} \bigl( \widehat{Z} ( \tau_m )
\bigr) - \phi_{\varepsilon} \bigl( \widehat{Z} ( 0 ) \bigr) & =& \int
_0^{\tau_m} \nabla\phi_{\varepsilon} \bigl(
\widehat{Z} ( s ) \bigr) \,\mathrm{d}\widehat{B} ( s )\nonumber
\\
&&{} + \sum_{j=1}^{n-1} \int
_0^{\tau_m} q_{\cdot,j} \cdot\nabla\phi
_{\varepsilon} \bigl( \widehat{Z} ( s ) \bigr) \,\mathrm{d} \widehat
{I}_j^{\infty,1}
( s )
\\
&&{} + \sum_{k,\ell= 1}^{n-1} \int
_0^{\tau_m} q^{ (2
)}_{\cdot, (k,\ell)} \cdot
\nabla\phi_{\varepsilon} \bigl( \widehat{Z} ( s ) \bigr) \,\mathrm{d}
\widehat{I}{}^{\infty,2}_{k,\ell} ( s ).\nonumber
\end{eqnarray}
Since $\widehat{B}$ has no drift and $\phi_{\varepsilon}$ and its
first derivatives
are bounded on each compact subset of $ ( \mathbb{R}_+
)^{n-1}$, the
definition of the stopping time $\tau_m$ implies that the stochastic
integral with respect to $\mathrm{d}\widehat{B}$ in~(\ref
{eqitophi}) has
zero expectation. To bound the other terms on the right-hand side
of~(\ref{eqitophi}) it is necessary to bound the directional
derivatives of $\phi_{\varepsilon}$; this is exactly what is done
in~\cite{RW}, pages~93--95. In particular, the results of~\cite{RW} give two
bounds. First, for every $j \in[n-1 ]$ there exists a
constant $\widehat{c}_j < \infty$ such that for all $x \in(
\mathbb{R}_+
)^{n-1}$ and all $\varepsilon\in( 0, 1 )$,
\[
q_{\cdot, j} \cdot\nabla\phi_{\varepsilon} ( x ) \geq-
\widehat{c}_j.
\]
Here the constant $\widehat{c}_j$ depends on $A$, $Q$, $\gamma$, and
$\delta
$, but does not depend on~$x$ nor~$\varepsilon$; see~\cite{RW}, equation~(24).
Next, for $j \in[ n - 1 ]$ let $\beta_j = \delta_j /
\llVert A^{-1} q_{\cdot, j} \rrVert$. Then for every $j \in[ n - 1
]$ there exists a constant $c_j > 0$ such that for all $x \in
( \mathbb{R}_+ )^{n-1}$ satisfying $\llVert x \rrVert
< \varepsilon\beta_j$,
%
%
\begin{equation}
\label{eqbdeps} q_{\cdot, j} \cdot\nabla\phi_{\varepsilon} ( x ) \geq-
c_j ( \ln\varepsilon+ 1 ).
\end{equation}
Here the constant $c_j$ depends on $A$, $Q$, $\gamma$, $\delta$ and
$\beta_j$, but does not depend on~$x$ nor~$\varepsilon$.
Note that for $\varepsilon$ small the term on the right-hand side
of~(\ref{eqbdeps}) is large and positive.
Furthermore, due to Assumption~\ref{mainass}(b) and the choice of
$\gamma$, the same arguments as in~\cite{RW}, pages~93--95, can be
repeated to bound the directional derivatives $q_{\cdot, (k,\ell
)}^{ (2 )} \cdot\nabla\phi_{\varepsilon}$. In
particular, for
every $ (k,\ell) \notin\mathcal{I}^{ (2 )}$ there
exists a constant $\widehat{c}_{ (k,\ell)} < \infty$
such that for
all $x\in( \mathbb{R}_+ )^{n-1}$ and all $\varepsilon
\in(0,1 )$,
\[
q_{\cdot, (k,\ell)}^{ (2 )} \cdot\nabla\phi_{\varepsilon} ( x ) \geq-
\widehat{c}_{ (k,\ell)}.
\]
Here the constant $\widehat{c}_{ (k,\ell)}$ depends on
$A, Q,
Q^{ (2 )}, \gamma$ and $\delta$, but does not depend on~$x$
nor $\varepsilon$. If $ (k,\ell) \in\mathcal{I}^{
(2 )}$,
then by definition $q_{\cdot, (k,\ell)}^{ (2
)}$ is
the zero vector,\vspace*{2pt} and thus $q_{\cdot, (k,\ell)}^{
(2
)} \cdot\nabla\phi_{\varepsilon} = 0$.
Plugging these bounds into~(\ref{eqitophi}) and taking expectation we
get that
%
%
\begin{eqnarray}
\label{eqphiest} &&\mathbb E \bigl[ \phi_{\varepsilon} \bigl( \widehat
{Z} ( \tau
_m ) \bigr) - \phi_{\varepsilon} \bigl( \widehat{Z} ( 0 ) \bigr)
\bigr]\nonumber
\\
&&\qquad \geq- ( \ln\varepsilon+ 1 ) \sum_{j=1}^{n-1}
c_j \mathbb E \biggl[ \int_0^{\tau_m}
\mathbf{1}_{ \{ \llVert
\widehat{Z} ( s
) \rrVert < \varepsilon\beta_j \}} \,\mathrm{d}\widehat{I}_j^{\infty,1} (
s ) \biggr]
\\
&&\quad\qquad{} - \sum_{j=1}^{n-1} \widehat{c}_j
\mathbb E \bigl[ \widehat{I}_j^{\infty,1} ( \tau
_m ) \bigr] - \sum_{ (k,\ell) \notin\mathcal{I}^{ (2 )}}
\widehat{c}_{ (k,\ell)} \mathbb E \bigl[ \widehat{I}{}^{\infty,2}_{k,\ell}
( \tau_m ) \bigr].\nonumber
\end{eqnarray}
The left-hand side of~(\ref{eqphiest}) is bounded as $\varepsilon
\to0$
since $\phi_{\varepsilon}$ is uniformly bounded on compact subsets of
$
(\mathbb{R}_+ )^{n-1}$, while the last two terms in~(\ref{eqphiest})
are independent of~$\varepsilon$.
So dividing~(\ref{eqphiest}) by $-( \ln\varepsilon+ 1
)$ and
letting $\varepsilon\to0$, we get that
\[
\lim_{\varepsilon\to0} \sum_{j=1}^{n-1}
c_j \mathbb E \biggl[ \int_0^{\tau_m}
\mathbf{1}_{ \{ \llVert \widehat{Z} ( s ) \rrVert < \varepsilon
\beta_j \}} \,\mathrm{d}\widehat{I}_j^{\infty,1} ( s
) \biggr] \leq0.
\]
Each term in the sum above is nonnegative and $c_j > 0$, so by Fatou's
lemma it follows that
\[
\int_0^{\tau_m} \mathbf{1}_{ \{ \widehat{Z}_{j'} ( s
) = 0,
j' \in[n-1 ] \} }
\,\mathrm{d}\widehat{I}_j^{\infty,1} ( s ) = 0
\]
for every $j \in[n-1 ]$ a.s. By letting $m \to\infty$ we
have that
\[
\int_0^{\infty} \mathbf{1}_{ \{ \widehat{Z}_{j'} ( s
) = 0,
j' \in[n-1 ] \} }
\,\mathrm{d}\widehat{I}_j^{\infty,1} ( s ) = 0
\]
for every $j \in[n-1 ]$ a.s. Finally, by using the backward
induction argument of~\cite{RW}, Lemma~5, it follows that with
probability one, for all $j \in[n - 1 ]$ and $J \subseteq
[n-1 ]$ such that $\llvert J \rrvert \geq2$ we have that
\[
\int_0^{\infty} \mathbf{1}_{ \{ \widehat{Z}_{j'} ( s
) = 0,
j' \in J \} }
\,\mathrm{d}\widehat{I}_j^{\infty,1} ( s ) = 0.
\]
Together with~(\ref{bdry1}), this implies~(\ref{eqbdryend1}).
\end{pf}

\subsection{General setup} \label{secgeneral}

In this last subsection, we introduce a much more general class of
particle systems which converge to appropriate multidimensional sticky
Brownian motions in the sense of Theorem \ref{mainthm}. We now allow
for nonexponential interarrival times between the jumps of the
particles and for dependence between the arrival times of the jumps for
different particles.

To define this more general class of particle systems, we introduce the
following parameters: $n \in\mathbb{N}$ for the number of particles as
before; $a>0$; $\lambda_i^L$ and $\lambda_i^R$ for $i\in[n]$;
$c_{i,i'}^{L,L}$, $c_{i,i'}^{L,R}$ and $c_{i,i'}^{R,R}$ for $i,i'\in
[n]$; and\vspace*{2pt} finally $\theta_{i,j}^L$ and $\theta_{i,j}^R$ for $i\in[n]$,
$j\in[n-1]$. We fix a value $M > 0$ of the scaling parameter. The
random variables and processes we define next all depend on $M$, but
for the sake of readability we mostly do not denote this dependence explicitly.

We let $ \{u^L(k), k\in\mathbb{N} \}$ and $ \{
u^R(k), k\in\mathbb{N}
\}$ be two independent sequences of i.i.d. random vectors with
values in $(0,\infty)^n$ (the interarrival times between jumps to the
left and to the right when there are no collisions), and for $i \in
[n ]$, $j \in[n-1 ]$, let $ \{ w_{i,j}^L
( k ), k\in\mathbb{N} \}$ and $ \{ w_{i,j}^R
(k ),
k\in\mathbb{N} \}$ be two\vspace*{2pt} independent families of i.i.d. random
variables taking values in $(0,\infty)$ (the interarrival times between
jumps to the left and to the right when there is a collision). We
assume that
\begin{eqnarray*}
\mathbb E\bigl[u_i^L(1)\bigr]&=& \biggl(a+
\frac{\lambda_i^L}{\sqrt{M}} \biggr)^{-1},\qquad\mathbb E\bigl[u_i^R(1)
\bigr]= \biggl(a+\frac{\lambda_i^R}{\sqrt{M}} \biggr)^{-1},
\\
\operatorname{cov}\bigl(u_i^L(1),u_{i'}^L(1)
\bigr)&=&c_{i,i'}^{L,L},\qquad \operatorname{cov}\bigl(u_i^R(1),u_{i'}^R(1)
\bigr)=c_{i,i'}^{R,R},
\\
\operatorname{cov}\bigl(u_i^L(1),u_{i'}^R(1)
\bigr)&=&c_{i,i'}^{L,R},
\\
\mathbb E\bigl[w_{i,j}^L(1)\bigr]&=& \bigl(
\theta^L_{i,j} \bigr)^{-1}\quad\mbox{and}\quad
\mathbb E\bigl[w_{i,j}^R(1)\bigr]= \bigl(
\theta^R_{i,j} \bigr)^{-1}.
\end{eqnarray*}

Next, define the corresponding partial sum processes
\begin{eqnarray*}
U_i^L (0 )&=&0,\qquad U_i^L (\ell)=
\sum_{k=1}^\ell u_i^L
(k ),
\\
U_i^R (0 )&=&0,\qquad U_i^R
(\ell)=\sum_{k=1}^\ell
u_i^R (k ),
\\
W_{i,j}^L (0 )&=&0,\qquad W_{i,j}^L (\ell)=
\sum_{k=1}^\ell w_{i,j}^L
(k ),
\\
W_{i,j}^R (0 )&=&0,\qquad W_{i,j}^R
(\ell)=\sum_{k=1}^\ell
w_{i,j}^R (k )
\end{eqnarray*}
for all $i\in[n]$, $j\in[n-1]$, and also the corresponding renewal processes
\begin{eqnarray*}
S_i^L (t )&=&\max\bigl\{k\geq0\dvtx  U_i^L
(k )\leq t \bigr\},\qquad S_i^R (t )=\max\bigl\{k\geq0\dvtx
U_i^R (k )\leq t \bigr\},
\\
T_{i,j}^L (t )&=&\max\bigl\{k\geq0\dvtx  W_{i,j}^L
(k )\leq t \bigr\},\qquad T_{i,j}^R (t )=\max\bigl\{k\geq0\dvtx
W_{i,j}^R (k )\leq t \bigr\}.
\end{eqnarray*}

We now define the particle system for any fixed value of the scaling
parameter $M>0$ according to
%
%
\begin{eqnarray}
\label{eqpartsysgen}
&& \mathrm{d} X_i^M ( t )\nonumber\hspace*{-10pt}
\\
&&\qquad = \frac{1}{\sqrt{M}} \mathbf{1}_{ \{X_k^M (t
)+({1}/{\sqrt{M}})<X_{k+1}^M (t ), k\in[n-1] \}} \,\mathrm{d}
\bigl(S_i^R (M t )-S_i^L (M t )
\bigr)\nonumber\hspace*{-10pt}
\\[-8pt]\hspace*{-10pt}
\\[-8pt]
&&\quad\qquad{} + \frac{1}{\sqrt{M}} \sum_{j=1}^{n-1}
\mathbf{1}_{ \{ X_i^M
( t ) + ({1}/{\sqrt{M}}) < X_{i+1}^M ( t ),
X_j^M ( t ) + ({1}/{\sqrt{M}}) = X_{j+1}^M ( t
) \}} \,\mathrm{d} T_{i,j}^R ( \sqrt{M} t)\nonumber\hspace*{-10pt}
\\
&&\quad\qquad{} - \frac{1}{\sqrt{M}} \sum_{j=1}^{n-1}
\mathbf{1}_{ \{ X_{i-1}^M
( t ) + ({1}/{\sqrt{M}}) < X_i^M ( t ), X_j^M
( t ) + ({1}/{\sqrt{M}}) = X_{j+1}^M ( t )
\}} \,\mathrm{d} T_{i,j}^L ( \sqrt{M} t),\nonumber\hspace*{-10pt}
\end{eqnarray}
for $i\in[n]$. Note that the particle configuration
\[
\bigl(X_1^M (t ),X_2^M (t ),
\ldots,X_n^M (t ) \bigr)
\]
is an element of the discrete wedge $\mathcal{W}^M$ for any $t\geq0$.

Intuitively, this general particle system behaves as follows. When
apart, the particles jump on the rescaled lattice $\mathbb{Z}/\sqrt
{M}$ with
jump rates of order $M$; these jumps are governed by the renewal
processes $S_i^L$ and $S_i^R$, $i \in[n ]$, and thus the
movements of the particles are not necessarily independent, and not
necessarily governed by Poisson processes. However, when a collision
occurs (i.e., two particles are on adjacent sites), then the system
experiences a slowdown, with the particles moving with jump rates of
order $\sqrt{M}$; these jumps are governed by the renewal processes
$T_{i,j}^L$ and $T_{i,j}^R$, $i \in[n ]$, $j \in[ n-1
]$, and thus the movements of the particles are independent, but
not necessarily governed by Poisson processes.

The particle system~(\ref{eqpartsysgen}) indeed generalizes~(\ref
{eqpartsys}), as the following parameter specifications show. If
$\lambda_i^L=\lambda_i^R=0$ for $i\in[n ]$,
$c_{i,i'}^{L,L}=c_{i,i'}^{L,R}=c_{i,i'}^{R,R}=0$ whenever $i\neq i'$,
$c_{i,i}^{L,L}=c_{i,i}^{R,R}=a^{-2}$ and $c_{i,i}^{L,R}=0$ for $i\in
[n]$, and all interarrival times above are independent exponential
random variables with appropriate means, then~(\ref{eqpartsysgen})
reduces to~(\ref{eqpartsys}).

For the extension of our convergence theorem to particle systems as
in~(\ref{eqpartsysgen}), we need the following moment assumption on
the interarrival times between jumps. This assumption is needed in
order to have uniform integrability of the appropriate sequences of
random variables.

%
%
\begin{asmp}\label{momass}
Assume that there exists $\delta>0$ such that
\begin{eqnarray*}
\sup_{M>0} \max_{i\in[n ]} \bigl(\mathbb E
\bigl[u_i^L (1 )^{2+\delta} \bigr]+\mathbb E
\bigl[u_i^R (1 )^{2+\delta
} \bigr] \bigr) &<& \infty,
\\
\sup_{M>0} \max_{i\in[n ], j\in[n-1]} \bigl(\mathbb E
\bigl[w^L_{i,j} (1 )^{2+\delta} \bigr]+\mathbb E
\bigl[w^R_{i,j} (1 )^{2+\delta} \bigr] \bigr) &<& \infty.
\end{eqnarray*}
\end{asmp}

Under Assumption~\ref{momass} we have the following convergence
result, which generalizes Theorem~\ref{mainthm}.

%
%
\begin{teo}\label{teogen}
Suppose that Assumptions~\ref{mainass} and~\ref{momass} hold, and
that the initial conditions $ \{ X^M (0 ), M>0 \}$
are deterministic and converge to a limit $x \in\mathcal{W}$ as $M
\to
\infty$. Then the laws of the paths of the particle systems $ \{
X^M (\cdot), M > 0 \}$ defined in~(\ref
{eqpartsysgen}) converge in $D ([0,\infty),\mathbb{R}^n
)$ to the
law of the unique weak solution of the system of SDEs
%
%
\begin{eqnarray}
\label{mainsdegen} \mathrm{d}X_i (t )& =&\mathbf{1}_{ \{X_1(t)<\cdots
<X_n(t) \}}
\bigl( \bigl(\lambda_i^R-\lambda_i^L
\bigr) \,\mathrm{d}t+a^{3/2} \,\mathrm{d}W_i (t ) \bigr)
\nonumber\\[-8pt]\\[-8pt]
&&{}+ \sum_{j=1}^{n-1} \mathbf{1}_{ \{X_j (t
)=X_{j+1}
(t ) \}}
v_{i,j} \,\mathrm{d}t\nonumber
\end{eqnarray}
for $i\in[n]$, taking values in $\mathcal{W}$ and starting from $x$.
Here, the vector $W= (W_1,W_2,\ldots,W_n )$ is a Brownian
motion in $\mathbb{R}^n$ with zero drift vector and diffusion matrix
given by
\[
\mathfrak{C}=(\mathfrak{c}_{i,i'})=\bigl
(c^{L,L}_{i,i'}+c^{L,R}_{i,i'}+c^{L,R}_{i',i}+c^{R,R}_{i,i'}
\bigr)
\]
and $v_{i,j}$ is as in~(\ref{eqvij}).
\end{teo}

The existence and uniqueness of a weak solution to the system of SDEs
given by~(\ref{mainsdegen}) is proven in Theorem~\ref{teouniq}, so
Theorem~\ref{teogen} is a consequence of Proposition~\ref
{propconvgen} below, which is the appropriate generalization of
Proposition~\ref{propconv} in Section~\ref{secconvbasic}.

As in Section~\ref{secconvbasic}, we need to study an appropriate
decomposition of the particle dynamics. For each $i\in[n ]$,
we write
\begin{eqnarray*}
X_i^M ( t ) & =& X_i^M ( 0 ) +
A_i^M ( t ) + \sum_{j=1}^{n-1}
C_{i,j}^{R,M} ( t ) - \sum_{j=1}^{n-1}
C_{i,j}^{L,M} ( t )
\\
&&{} + \sum_{j=1}^{n-1} \Delta_{i,j}^{R,M}
( t ) - \sum_{j=1}^{n-1}
\Delta_{i,j}^{L,M} ( t ),
\end{eqnarray*}
where now
\begin{eqnarray*}
A_i^M ( t )
&:=&\frac{1}{\sqrt{M}} \int_0^t
\mathbf{1}_{ \{X_k^M
(s
)+({1}/{\sqrt{M}})<X_{k+1}^M (s ), k\in[n-1] \}
} \,\mathrm{d} \bigl(S_i^R (M s
)-S_i^L (M s ) \bigr),
\\
C_{i,j}^{R,M} ( t )&:=&\theta_{i,j}^R
I_{i,j}^{R,M}
\\
&:=& \theta_{i,j}^R \int
_0^t \mathbf{1}_{ \{X_i^M (s
)+({1}/{\sqrt{M}})<X_{i+1}^M (s ),X_j^M (s )+
({1}/{\sqrt
{M}})=X_{j+1}^M (s ) \}} \,\mathrm{d}s,
\\
\Delta_{i,j}^{R,M} (t )
&:=&\frac{1}{\sqrt{M}} \int_0^t
\mathbf{1}_{ \{X_i^M (s )+({1}/{\sqrt
{M}})<X_{i+1}^M (s ),X_j^M (s )+({1}/{\sqrt
{M}})=X_{j+1}^M (s ) \}} \,\mathrm{d}\overline{T}_{i,j}^R (
\sqrt{M} s ),
\\
\overline{T}_{i,j}^R (t )&:=&T_{i,j}^R
(t )-\theta_{i,j}^R t,
\end{eqnarray*}
and the processes $C_{i,j}^{L,M}$, $I^{L,M}_{i,j}$, $\Delta
_{i,j}^{L,M}$ and $\overline{T}_{i,j}^L$ are defined similarly to
$C_{i,j}^{R,M}$, $I_{i,j}^{R,M}$, $\Delta_{i,j}^{R,M}$ and $\overline
{T}_{i,j}^R$, respectively. The following proposition is the
appropriate generalization of Proposition~\ref{propconv} to the
present framework.

%
%
\begin{prop}\label{propconvgen}
Suppose that Assumptions~\ref{mainass} and~\ref{momass} hold, and
that the initial conditions $ \{ X^M (0 ), M>0 \}$
are deterministic and converge to a limit $x\in\mathcal{W}$ as
$M\rightarrow\infty$. Then the family
%
%
\begin{equation}
\label{eqfamgen} \bigl\{ \bigl(X^M,A^M,I^{L,M},I^{R,M},
\Delta^{L,M},\Delta^{R,M} \bigr), M>0 \bigr\}
\end{equation}
is tight in $D^{4n^2-2n}$. Moreover, every limit point
\[
\bigl(X^\infty,A^\infty,I^{L,\infty},I^{R,\infty},\Delta
^{L,\infty
},\Delta^{R,\infty} \bigr)
\]
satisfies the following for each $i \in[n ]$:
%
%
\begin{eqnarray}
 X_i^\infty( \cdot) & =& \int
_0^\cdot\mathbf{1}_{ \{
X_1^\infty(s )<\cdots<X_n^{\infty} ( s )
\}} \bigl( \bigl(
\lambda_i^R - \lambda_i^L
\bigr) \,\mathrm{d}s + a^{3/2} \,\mathrm{d}W_i (s ) \bigr)
\nonumber\\[-8pt]\label{eqconvXgen} \\[-8pt]
&&{} +\sum_{j=1}^{n-1} v_{i,j} \int
_0^\cdot\mathbf{1}_{ \{
X_j^{\infty} ( s ) = X_{j+1}^\infty( s )
\}
} \,\mathrm{d}s,\nonumber
\\
A_i^\infty(\cdot) & =& \int_0^\cdot
\mathbf{1}_{ \{
X_1^\infty(s )<\cdots<X_n^{\infty} (s )
\}} \bigl( \bigl(\lambda_i^R-
\lambda_i^L \bigr) \,\mathrm{d}s+a^{3/2} \,\mathrm
{d}W_i (s ) \bigr), \label{eqconvAgen}
\\
I_{i,j}^{L,\infty} (\cdot) & =& \int_0^\cdot
\mathbf{1}_{ \{X_j^{\infty} (s )=X_{j+1}^\infty
(s )
\}} \,\mathrm{d}s, \qquad j \in[n-1 ] \setminus\{ i-1 \},
\label{eqconvCLgen}
\\
I_{i,j}^{R,\infty} (\cdot) & =& \int_0^\cdot
\mathbf{1}_{ \{X_j^{\infty} (s )=X_{j+1}^\infty
(s )
\}} \,\mathrm{d}s, \qquad j \in[n-1 ] \setminus\{ i \},
\label{eqconvCRgen}
\\
I_{i,i-1}^{L,\infty} ( \cdot) & =& I_{i,i}^{R,\infty}
( \cdot) = 0,
\nonumber
\\
\Delta^{L,\infty}_{i,j}(\cdot) & =& \Delta^{R,\infty}_{i,j}(
\cdot)=0, \qquad j \in[n-1 ],
\nonumber
\end{eqnarray}
with a Brownian motion $W= (W_1,W_2,\ldots,W_n )$ as in the
statement of Theorem~\ref{teogen}.
\end{prop}

\begin{pf}
One can proceed as in the proof of Proposition~\ref{propconv}, so we
only explain the arguments which are different. First, note that
Theorem~14.6 in~\cite{Bi} and its proof extend to the case of the
multidimensional renewal processes
\[
\bigl\{ S_i^L \bigr\}_{i \in[n ]}, \qquad\bigl\{
S_i^R \bigr\}_{i \in[n ]}, \qquad\bigl\{
T_{i,j}^L \bigr\}_{i
\in[n ], j \in[n-1 ]}, \qquad\bigl\{
T_{i,j}^R \bigr\}_{i \in[n ], j \in[n-1 ]},
\]
yielding the joint convergence of
\begin{eqnarray*}
&\displaystyle \bigl\{ \bigl( M^{-1/2}\bigl(S_i^R(M
t)-S_i^L(M t)\bigr), t \geq0 \bigr) \bigr
\}_{i \in[n ]},&
\\
& \displaystyle\bigl\{ \bigl( M^{-1/4} \overline{T}_{i,j}^L(
\sqrt{M} t), t \geq0 \bigr) \bigr\}_{i\in[n ], j \in[n -1 ]},&
\end{eqnarray*}
and
\[
\bigl\{ \bigl( M^{-1/4} \overline{T}_{i,j}^R(
\sqrt{M} t), t \geq0 \bigr) \bigr\}_{i\in[n ], j \in[n -1 ]}
\]
to appropriate Brownian motions. The rest of steps~1~and~2 in the
proof of Proposition~\ref{propconv} carry over to the present setting
in a straightforward manner.

Now, one needs to show that every limit point
\[
\bigl(X^\infty,A^\infty,I^{L,\infty},I^{R,\infty},\Delta
^{L,\infty
},\Delta^{R,\infty} \bigr)
\]
satisfies
\[
\bigl\langle X^\infty_i,X^\infty_{i'}
\bigr\rangle(\cdot)= \bigl\langle A^\infty_i,A^\infty_{i'}
\bigr\rangle(\cdot) =a^3 \mathfrak{c_{i,i'}} \int
_0^\cdot\mathbf{1}_{ \{
X_1^\infty
(s )<\cdots<X_n^{\infty} (s ) \}} \,\mathrm{d}s.
\]
To this end, one can first proceed as in step~3 in the proof of
Proposition \ref{propconv} to show that
%
%
\begin{equation}
\label{qv1} \mathrm{d} \bigl\langle X^\infty_i,X^\infty_{i'}
\bigr\rangle= \frac
{\mathfrak{c}_{i,i'}}{\mathfrak{c}_{j,j}} \,\mathrm{d} \bigl\langle
X^\infty_j
\bigr\rangle, \qquad i,i'\in[n], j\in[n].
\end{equation}
Next, one can invoke the Portmanteau theorem as before to conclude that
for all $i\in[n]$ and $0 \leq t_1 < t_2$,
%
%
\begin{equation}
\label{qv2} 
a^3 \mathfrak{c_{i,i}} \int
_{t_1}^{t_2} \mathbf{1}_{ \{
X_1^\infty(s )<\cdots<X_n^{\infty} (s
) \}} \,\mathrm{d}s
\leq\bigl\langle X^\infty_i \bigr\rangle(t_2) -
\bigl\langle X^\infty_i \bigr\rangle(t_1)
\leq a^3 \mathfrak{c_{i,i}} (t_2-t_1).\hspace*{-20pt}
\end{equation}
Moreover, since the measures $\mathrm{d} \langle X^\infty
_j
\rangle$, $j \in[n-1 ]$, assign zero mass to the sets $\{
t\geq0\dvtx  X^\infty_j(t)=X^\infty_{j+1}(t)\}$, $j\in[n-1]$, respectively,
(\ref{qv1}) and (\ref{qv2}) suffice to identify all quadratic
covariation processes $ \langle X^\infty_i,X^\infty_{i'}
\rangle$, $i,i'\in[n]$. Similarly, one can identify the bounded
variation parts of the processes $A^\infty_i$, $i\in[n]$, as multiples
of the quadratic variation processes $ \langle X^\infty_i
\rangle$, $i\in[n]$, respectively. The rest of the proof can be carried
out by following the arguments in step~4 in the proof of
Proposition~\ref{propconv}.
\end{pf}


\section*{Acknowledgments}

The authors thank Soumik Pal for many helpful discussions in the course
of the preparation of this paper,
and two anonymous referees for numerous useful suggestions that have
helped catch mistakes and improve the exposition of the paper.


akaldyti doi

%

\printaddresses

\end{document}